\documentclass{amsart}
\usepackage{amssymb}
\usepackage{amsmath,amsthm}
\usepackage{afterpage}
\usepackage{graphicx}
\title[Exotic rational elliptic surfaces without 1-handles]{Exotic rational elliptic surfaces\\without 1-handles}
\author[Kouichi Yasui]{Kouichi Yasui}
\date{May 9, 2007}
\address{Department~of~Mathematics, Graduate~School~of~Science, Osaka~University, Toyonaka, Osaka 560-0043, Japan}
\email{kyasui@cr.math.sci.osaka-u.ac.jp}
\thanks{The author is partially supported by JSPS Research Fellowships for Young Scientists.}

\newtheorem{theorem}{Theorem}[section]
\newtheorem{proposition}[theorem]{Proposition}
\newtheorem{lemma}[theorem]{Lemma}
\newtheorem{corollary}[theorem]{Corollary}
\theoremstyle{definition}
\newtheorem{definition}[theorem]{Definition}
\newtheorem{remark}[theorem]{Remark}
\newtheorem{problem}[theorem]{Problem}
\newtheorem{ack}{Acknowledgement}

\begin{document}

\begin{abstract}
Harer, Kas and Kirby have conjectured that every handle decomposition of the elliptic surface $E(1)_{2,3}$ requires both $1$- and $3$-handles. 
In this article, we construct a smooth $4$-manifold which has the same Seiberg-Witten invariant as $E(1)_{2,3}$ and admits neither $1$- nor $3$-handles, by using rational blow-downs and Kirby calculus. Our manifold gives the first example of either a counterexample to the Harer-Kas-Kirby conjecture or a homeomorphic but non-diffeomorphic pair of simply connected closed smooth $4$-manifolds with the same non-vanishing Seiberg-Witten invariants.
\end{abstract}

\maketitle

\section{Introduction}
It is a basic problem in $4$-dimensional topology to classify smooth structures on $4$-manifolds.
Constructions of exotic smooth structures on $4$-manifolds with small Euler characteristics are currently in rapid progress (see, for example, Park \cite{P2}, Stipsicz-Szab\'{o} \cite{SS}, Fintushel-Stern \cite{FS2}, Park-Stipsicz-Szab\'{o} \cite{PSS} and Akhmedov-Park \cite{AP}). However, it is still unknown whether or not $\mathbf{S}^4$ and $\mathbf{CP}^2$ admit an exotic smooth structure. If such a structure exists, then each handle decomposition of it has at least either a $1$- or $3$-handle (see Proposition~\ref{S^4}). To the contrary, many classical simply connected closed smooth $4$-manifolds are known to 
admit neither $1$- nor $3$-handles in their handle decompositions (cf.~Gompf-Stipsicz \cite{GS}).
Problem 4.18 in Kirby's problem list \cite{Ki} is the following: ``Does every simply connected, closed $4$-manifold 
have a handlebody decomposition without $1$-handles? Without $1$- and $3$-handles?'' 
The elliptic surfaces $E(n)_{p,q}$ are candidates of counterexamples to Problem 4.18. 
It is not known whether or not the simply connected closed smooth $4$-manifold $E(n)_{p,q}\,
(n\colon \text{arbitrary}, p,q\geq 2,\, \gcd (p,q)=1)$ 
admits a handle decomposition without $1$-handles (cf.~Gompf \cite{G} and Gompf-Stipsicz \cite{GS}). 
In particular, Harer, Kas and Kirby have conjectured in \cite{HKK} that 
every handle decomposition of $E(1)_{2,3}$ requires at least a $1$-handle. 
Note that by considering dual handle decompositions, their conjecture is equivalent to the assertion that 
$E(1)_{2,3}$ requires both $1$- and $3$-handles.

In this article we construct the following smooth $4$-manifolds 
by using rational blow-downs and Kirby calculus. 
\begin{theorem}
{\normalfont{(1)}} For $q=3,5$, there exists a smooth $4$-manifold $E_q$ with the following properties:\\
{\normalfont{(a)}} $E_q$ is homeomorphic to $E(1)_{2,q}$;\\
{\normalfont{(b)}} $E_q$ has the same Seiberg-Witten invariant as $E(1)_{2,q}$;\\
{\normalfont{(c)}} $E_q$ admits a handle decomposition without $1$-handles, namely,
\begin{equation*}
E_q=\text{one $0$-handle} \cup \text{twelve $2$-handles} \cup \text{two $3$-handles} \cup \text{one $4$-handle}.
\end{equation*}
{\normalfont{(2)}} There exists a smooth $4$-manifold $E'_3$ with the following properties:\\
{\normalfont{(a)}} $E'_3$ is homeomorphic to $E(1)_{2,3}$;\\
{\normalfont{(b)}} $E'_3$ has the same Seiberg-Witten invariant as $E(1)_{2,3}$;\\
{\normalfont{(c)}} $E'_3$ admits a handle decomposition without $1$- and $3$-handles, namely,
\begin{equation*}
E'_3=\text{one $0$-handle} \cup \text{ten $2$-handles} \cup \text{one $4$-handle}.
\end{equation*}
\end{theorem}

As far as the author knows, $E_q$ and $E'_3$ are the first examples in the following sense: 
If $E_q$ (resp.~$E'_3$) is diffeomorphic to $E(1)_{2,q}$ (resp.~$E(1)_{2,3}$), then the above handle decomposition of 
$E(1)_{2,q}$ ($=E_q$ [resp.~$E'_3$]) is the first example 
which has no $1$-handles. 
Otherwise, i.e., if $E_q$ (resp.~$E'_3$) is not diffeomorphic to $E(1)_{2,q}$ (resp.~$E(1)_{2,3}$), then $E_q$ (resp.~$E'_3$) and $E(1)_{2,q}$ (resp.~$E(1)_{2,3}$) are the first homeomorphic but non-diffeomorphic examples which are 
simply connected closed smooth $4$-manifolds with the same non-vanishing Seiberg-Witten invariants.

An affirmative solution to the Harer-Kas-Kirby conjecture 
implies that both $E_3$ and $E'_3$ are not diffeomorphic to $E(1)_{2,3}$, 
though these three have the same Seiberg-Witten invariants. 
In this case, the minimal number of $1$-handles in handle decompositions does detect the difference of their smooth structures.

Our construction is inspired by rational blow-down constructions of exotic smooth structures on $\mathbf{CP}^2\# n\overline{\mathbf{C}\mathbf{P}^2}\,(5\leq n\leq 8)$ by Park \cite{P2}, Stipsicz-Szab\'{o} \cite{SS}, Fintushel-Stern \cite{FS2} and Park-Stipsicz-Szab\'{o} \cite{PSS}. Our method is different from theirs since, firstly, we use Kirby calculus to perform rational blow-downs, whereas they used elliptic fibrations on $E(1)$ (and knot surgeries), secondly, they did not examine handle decompositions. 
\begin{ack} \normalfont 
The author wishes to express his deeply gratitude to his adviser, 
Professor Hisaaki Endo, for encouragement and many useful suggestions. 
He would like to thank Professors Selman Akbulut, Kazunori Kikuchi, Ronald J. Stern and Yuichi Yamada for helpful comments and discussions. 
Kikuchi's theorem \cite[Theorem 4]{K} partially gave him the idea of the construction. Yamada gave him interesting questions (cf.~Remark~\ref{Yamada}). 
\end{ack}
\section{Rational blow-down}
In this section we review the rational blow-down introduced 
by Fintushel-Stern \cite{FS1}. For details, see also Gompf-Stipsicz \cite{GS}.

\begin{figure}[htbp]
\begin{center}
\includegraphics[width=4.5in]{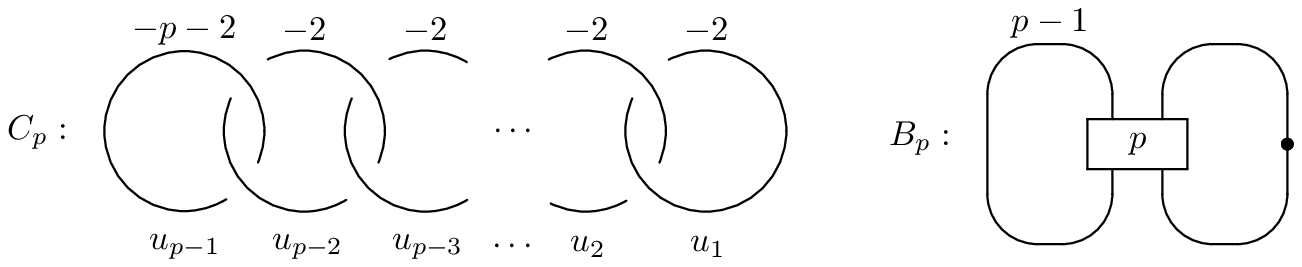}
\caption{}
\label{C_p}
\end{center}
\end{figure}

Let $C_p$ and $B_p$ be the smooth $4$-manifolds defined by Kirby diagrams in Figure~\ref{C_p}, and $u_1,\dots,u_{p-1}$ elements of $H_2(C_p;\mathbf{Z})$ given by corresponding $2$-handles in the figure such that $u_i\cdot u_{i+1}=+1$ $(1\leq i\leq p-2)$.
The boundary $\partial C_p$ of $C_p$ is diffeomorphic to the lens space $L(p^2,1-p)$ and to the boundary $\partial B_p$ of $B_p$. 
The following lemma is well known. 
\begin{lemma}
$(1)$ $\pi_1(C_p)=0$, $\pi_1(B_p)=\mathbf{Z}_p$ and $\pi_1(L(p^2,1-p))=\mathbf{Z}_{p^2}$.\\
$(2)$ $H_2(C_p;\mathbf{Z})=\oplus_{p-1} \mathbf{Z}$ and $H_2(B_p;\mathbf{Z})=H_2(L(p^2,1-p);\mathbf{Z})=0$
\end{lemma}
Suppose that $C_p$ embeds in a smooth $4$-manifold $X$. 
The smooth $4$-manifold $X_{(p)}:=(X-\text{int }C_p)\cup _{L(p^2,1-p)}B_p$ is called the rational blow-down of $X$ along $C_p$. Note that $X_{(p)}$ is uniquely determined up to diffeomorphism by a fixed pair $(X,C_p)$. 
This operation preserves $b_2^+$, decreases $b_2^-$, may create torsions in the first homology group, and has the following relation with the logarithmic transformation.
\begin{theorem}[{Fintushel-Stern \cite{FS1}, cf.~Gompf-Stipsicz \cite{GS}}]\label{th:2.1}
Suppose that a smooth $4$-manifold $X$ contains a cusp neighborhood, that is, a $0$-handle with a $2$-handle attached along a $0$-framed right trefoil knot. 
Let $X_{p}$ be the smooth $4$-manifold obtained from $X$ by 
performing a logarithmic transformation of multiplicity $p$ in the cusp neighborhood. 
Then there exists a copy of $C_p$ in 
$X\# (p-1)\overline{\mathbf{C}\mathbf{P}^2}$ such that the rational blow-down of 
$X\# (p-1)\overline{\mathbf{C}\mathbf{P}^2}$ along the copy of $C_p$ is diffeomorphic to $X_p$.
\end{theorem}
Let $E(n)$ be the simply connected elliptic surface with Euler characteristic $12n$ 
and with no multiple fibers, and $E(n)_{p_1,\dots,p_k}$ the elliptic surface obtained from $E(n)$ by performing logarithmic transformations of multiplicities $p_1,\dots,p_k$. We denote $h,e_1,e_2,\dots,e_n$ as a canonical orthogonal basis of $H_2(\mathbf{CP}^2\# n\overline{\mathbf{C}\mathbf{P}^2};\mathbf{Z})=H_2(\mathbf{CP}^2;\mathbf{Z})\oplus _n H_2(\overline{\mathbf{C}\mathbf{P}^2};\mathbf{Z})$ such that $h^2=1$ and $e_1^2=e_2^2=\dots=e_n^2=-1$. 

Since there is a diffeomorphism $E(1)_p\to E(1)=\mathbf{CP}^2\# 9\overline{\mathbf{C}\mathbf{P}^2}$ which maps the class of a regular fiber of $E(1)_p$ to $p(3h-e_1-e_2-\dots -e_9)\in H_2(\mathbf{CP}^2\# 9\overline{\mathbf{C}\mathbf{P}^2};\mathbf{Z})$ (cf.~Etg\"u-Park \cite[page 680]{EP}, Gompf-Stipsicz \cite{GS}), Theorem~\ref{th:2.1} gives us the following corollary. 
\begin{corollary}\label{cor:2.2}
For each natural number $p$ and $q$, the elliptic surface $E(1)_{p,q}$ is obtained from $\mathbf{CP}^2\# (8+q)\overline{\mathbf{C}\mathbf{P}^2}$ by rationally blowing down along a certain copy $_pC_q$ of $C_q$ such that $u_{1},\dots,u_{q-1}$ satisfy 
\begin{align*}
&u_{1}=e_{7+q}-e_{8+q},\,u_{2}=e_{6+q}-e_{7+q},\,\dots,\,u_{q-2}=e_{10}-e_{11},\\
&u_{q-1}=p(3h-e_1-e_2-\dots -e_9)-2e_{10}-e_{11}-e_{12}-\dots-e_{8+q}
\end{align*}
as elements of $H_2(\mathbf{CP}^2\# (8+q)\overline{\mathbf{C}\mathbf{P}^2};\mathbf{Z})$.
\end{corollary}
\begin{remark}
$E(1)_{p,q}$ is homeomorphic but non-diffeomorphic to $E(1)$, in the case $p,q\geq 2$ and $\gcd (p,q)=1$ (cf.~Gompf-Stipsicz \cite{GS}).
\end{remark}
\section{Construction}
In this section we construct $E_3$, $E_5$ and $E_3'$, and prove Theorem 1.1.(1)(a)(c) and (2)(a)(c). In Kirby diagrams, we write the second homology classes given by $2$-handles, instead of usual framings. Note that the square of the homology class given by a 2-handle is equal to the usual framing. We do not draw (whole) Kirby diagrams of $E_3,E_5,E'_3$ and the other manifolds appeared in the following construction. However, one can easily draw whole diagrams. 

\begin{figure}[p]
\begin{center}
\begin{minipage}{.45\linewidth}
\begin{center}
\includegraphics[scale=0.77]{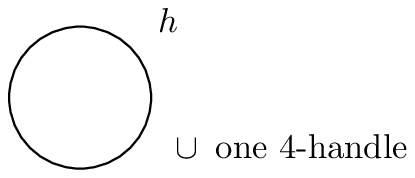}
\caption{$\mathbf{CP}^2$}
\label{4.1}
\end{center}
\end{minipage}
\begin{minipage}{.45\linewidth}
\begin{center}
\includegraphics[scale=0.77]{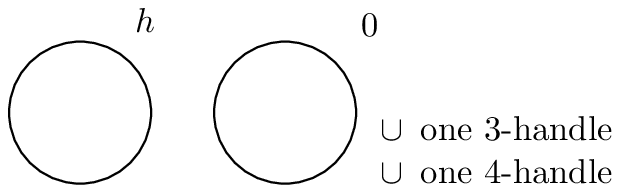}
\caption{$\mathbf{CP}^2$}
\label{4.2}
\end{center}
\end{minipage}
\end{center}
\end{figure}

\begin{figure}[htbp]
\begin{center}
\begin{minipage}{.45\linewidth}
\begin{center}
\includegraphics[scale=0.77]{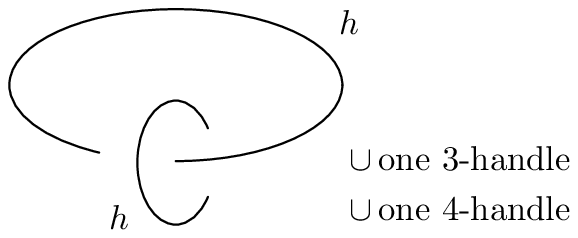}
\caption{$\mathbf{CP}^2$}
\label{4.3}
\end{center}
\end{minipage}
\begin{minipage}{.45\linewidth}
\begin{center}
\includegraphics[scale=0.77]{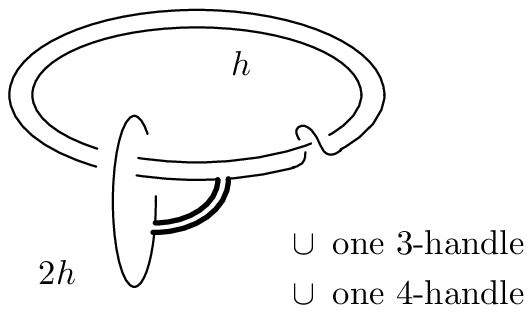}
\caption{$\mathbf{CP}^2$}
\label{4.4}
\end{center}
\end{minipage}
\end{center}
\end{figure}

\begin{figure}[htbp]
\begin{center}
\begin{minipage}{.45\linewidth}
\begin{center}
\includegraphics[scale=0.77]{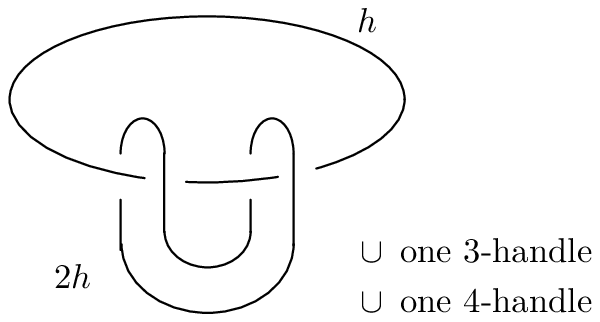}
\caption{$\mathbf{CP}^2$}
\label{4.5}
\end{center}
\end{minipage}
\begin{minipage}{.45\linewidth}
\begin{center}
\includegraphics[scale=0.77]{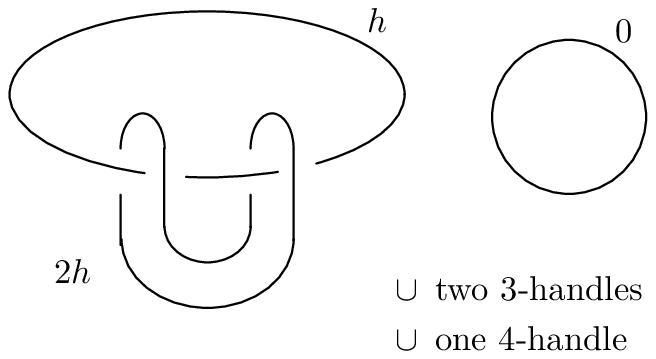}
\caption{$\mathbf{CP}^2$}
\label{4.6}
\end{center}
\end{minipage}
\end{center}
\end{figure}

\begin{figure}[htbp]
\begin{center}
\begin{minipage}{.45\linewidth}
\begin{center}
\includegraphics[scale=0.78]{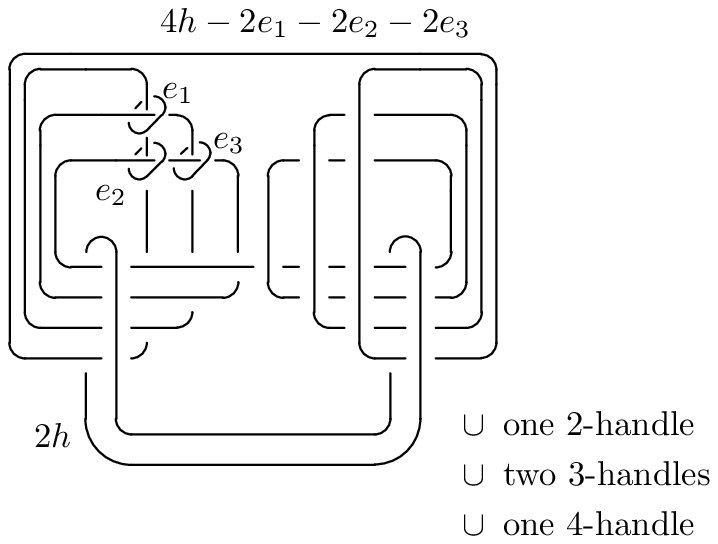}
\caption{$\mathbf{CP}^2\# 3\overline{\mathbf{C}\mathbf{P}^2}$}
\label{4.11}
\end{center}
\end{minipage}
\begin{minipage}{.45\linewidth}
\begin{center}
\includegraphics[scale=0.78]{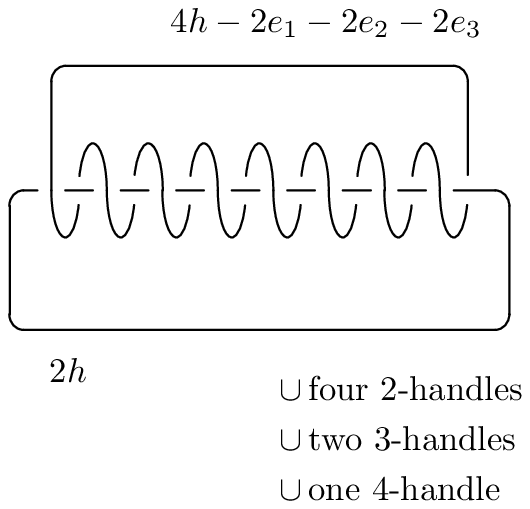}
\caption{$\mathbf{CP}^2\# 3\overline{\mathbf{C}\mathbf{P}^2}$}
\label{4.12}
\end{center}
\end{minipage}
\end{center}
\end{figure}

\begin{figure}[htbp]
\begin{center}
\begin{minipage}{.45\linewidth}
\begin{center}
\includegraphics[scale=0.78]{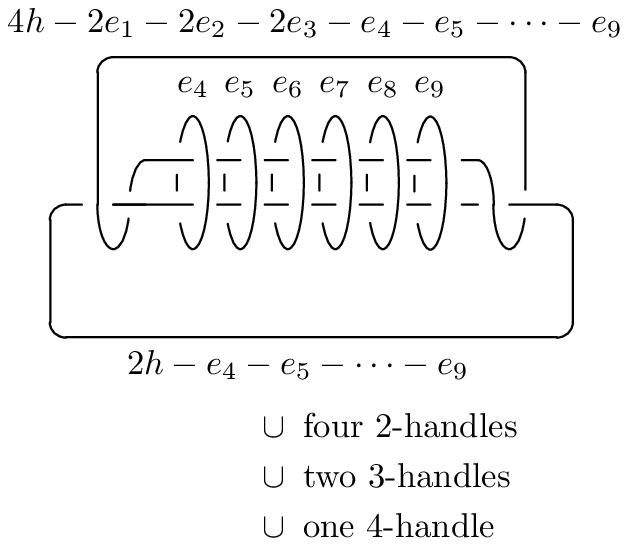}
\caption{$\mathbf{CP}^2\# 9\overline{\mathbf{C}\mathbf{P}^2}$}
\label{4.13}
\end{center}
\end{minipage}
\begin{minipage}{.45\linewidth}
\begin{center}
\includegraphics[scale=0.78]{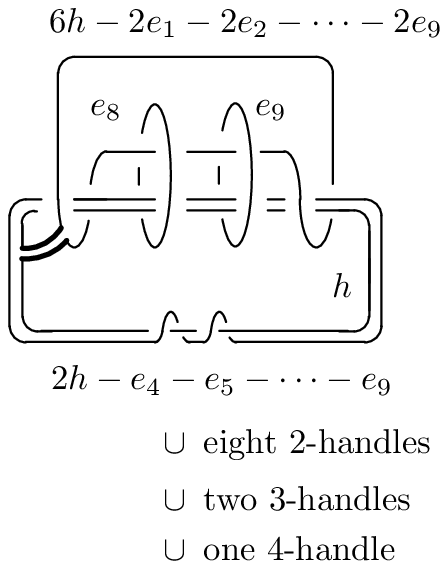}
\caption{$\mathbf{CP}^2\# 9\overline{\mathbf{C}\mathbf{P}^2}$}
\label{4.14}
\end{center}
\end{minipage}
\end{center}
\end{figure}

\begin{figure}[htbp]
\begin{center}
\includegraphics[scale=0.81]{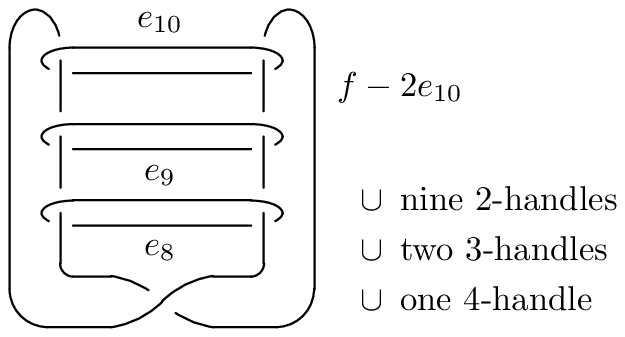}
\caption{$\mathbf{CP}^2\# 10\overline{\mathbf{C}\mathbf{P}^2}$}
\label{5.2}
\end{center}
\end{figure}

\begin{figure}[htbp]
\begin{center}
\includegraphics[scale=0.91]{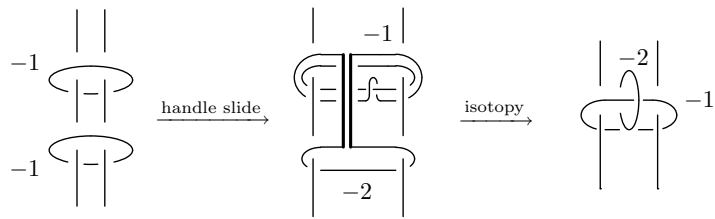}
\caption{Handle slide}
\label{ex-handleslide}
\end{center}
\end{figure}

\begin{figure}[htbp]
\begin{center}
\begin{minipage}{.45\linewidth}
\begin{center}
\includegraphics[scale=0.75]{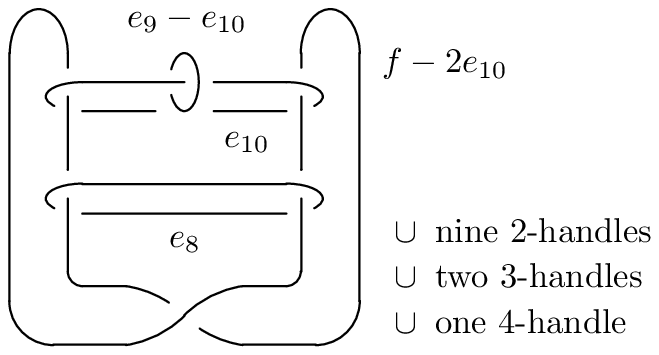}
\caption{$\mathbf{CP}^2\# 10\overline{\mathbf{C}\mathbf{P}^2}$}
\label{5.3}
\end{center}
\end{minipage}
\begin{minipage}{.45\linewidth}
\begin{center}
\includegraphics[scale=0.75]{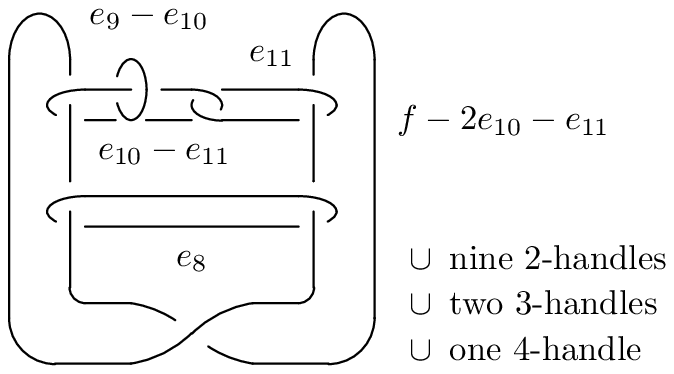}
\caption{$\mathbf{CP}^2\# 11\overline{\mathbf{C}\mathbf{P}^2}$}
\label{5.4}
\end{center}
\end{minipage}
\end{center}
\end{figure}

We begin with a construction of a cusp neighborhood in $\mathbf{CP}^2\# 9\overline{\mathbf{C}\mathbf{P}^2}$ such that its embedding into $\mathbf{CP}^2\# 9\overline{\mathbf{C}\mathbf{P}^2}$ has the same homological properties as that of the regular neighborhood of a cusp fiber of $E(1)_2$. We do not know if these embeddings into $\mathbf{CP}^2\# 9\overline{\mathbf{C}\mathbf{P}^2}$ are the same up to diffeomorphism.
\begin{lemma}\label{lem:3.1}
$\mathbf{CP}^2\# 9\overline{\mathbf{C}\mathbf{P}^2}$ admits the handle decomposition drawn in Figure~\ref{cusp}. Here $f$ denotes $6h-2e_1-2e_2-\dots-2e_9\in H_2(\mathbf{CP}^2\# 9\overline{\mathbf{C}\mathbf{P}^2};\mathbf{Z})$.
\begin{figure}[ht!]
\begin{center}
\includegraphics[width=2.0in]{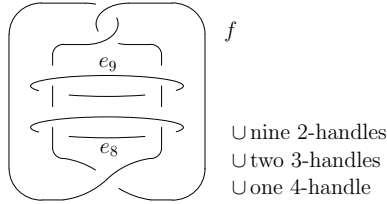}
\end{center}
\caption{$\mathbf{CP}^2\# 9\overline{\mathbf{C}\mathbf{P}^2}$}
\label{cusp}
\end{figure}
\end{lemma}
\begin{proof}
We firstly create two $2$-handles with framings $2h$ and $4h$ in a Kirby diagram of $\mathbf{CP}^2$. Figure~\ref{4.1} is a basic Kirby diagram of $\mathbf{CP}^2$. Introducing a $2$-handle/$3$-handle pair gives Figure~\ref{4.2}. Handle slides and isotopies yield Figure~\ref{4.5}\,(Pairs of bold lines in figures denote `bands'):
\begin{equation*}
\text{Figure~\ref{4.2}\: $\xrightarrow{0+h}$\: Figure~\ref{4.3}\:$\xrightarrow{h+h}$\: Figure~\ref{4.4}\: $\xrightarrow{\text{isotopy}}$\: Figure~\ref{4.5}.}
\end{equation*}
Creating a $2$-handle/$3$-handle pair gives Figure~\ref{4.6}. Handle slides produce Figure~\ref{4.10}:
\begin{equation*}
\text{Figure~\ref{4.6}\: $\xrightarrow{0+h}$\: Figure~\ref{4.7}\:$\xrightarrow{h+h}$\: Figure~\ref{4.8}\: $\xrightarrow{2h+h}$\: Figure~\ref{4.9}\: $\xrightarrow{3h+h}$\: Figure~\ref{4.10}.}
\end{equation*}
We secondly blow up $\mathbf{CP}^2$ nine times: 
\begin{equation*}
\text{Figure~\ref{4.10}\: $\xrightarrow{\text{three blow-ups}}$\: Figure~\ref{4.11}\:$\xrightarrow{\text{isotopy}}$\: Figure~\ref{4.12}\: $\xrightarrow{\text{six blow-ups}}$\: Figure~\ref{4.13}.}
\end{equation*}
We lastly make a handle addition $(4h-2e_1-2e_2-2e_3-e_4-e_5-\dots-e_{9})+(2h-e_4-e_5-\dots-e_9)$. This leads to Figure~\ref{4.14}, and an isotopy gives Figure~\ref{cusp}.
\end{proof}
\begin{proposition}\label{prop:3.2}
\text{\normalfont{(1)}}\, $\mathbf{CP}^2\# 11\overline{\mathbf{C}\mathbf{P}^2}$ 
admits the handle decomposition drawn in Figure~\ref{5.1}. In particular $\mathbf{CP}^2\# 11\overline{\mathbf{C}\mathbf{P}^2}$ contains the copy of $C_3$ drawn in the figure. The elements $u_1,u_2\in H_2(\mathbf{CP}^2\# 11\overline{\mathbf{C}\mathbf{P}^2};\mathbf{Z})$ given by this copy of $C_3$ are the same as that given by $_2C_3$.
\begin{figure}[ht!]
\begin{center}
\includegraphics[scale=0.75]{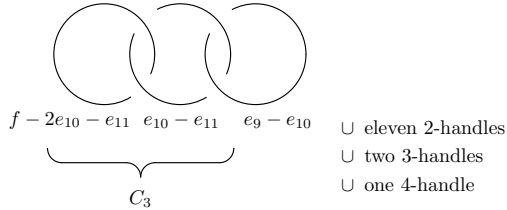}
\caption{$\mathbf{CP}^2\# 11\overline{\mathbf{C}\mathbf{P}^2}$}
\label{5.1}
\end{center}
\end{figure}\\
\text{\normalfont{(2)}}\, $\mathbf{CP}^2\# 13\overline{\mathbf{C}\mathbf{P}^2}$ admits the handle decomposition drawn in Figure~\ref{6.1}. In particular $\mathbf{CP}^2\# 13\overline{\mathbf{C}\mathbf{P}^2}$ contains the copy of $C_5$ drawn in the figure. The elements $u_1,\dots,u_4\in H_2(\mathbf{CP}^2\# 13\overline{\mathbf{C}\mathbf{P}^2};\mathbf{Z})$ given by this copy of $C_5$ are the same as that given by $_2C_5$.
\begin{figure}[ht!]
\begin{center}
\includegraphics[scale=0.70]{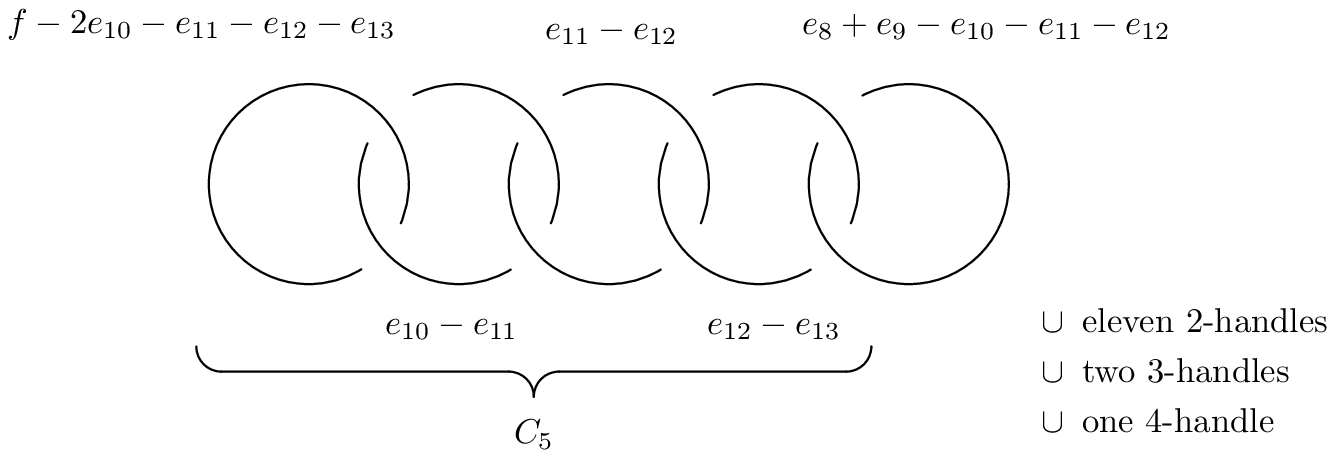}
\caption{$\mathbf{CP}^2\# 13\overline{\mathbf{C}\mathbf{P}^2}$}
\label{6.1}
\end{center}
\end{figure}
\end{proposition}
\proof
Firstly we give a proof for (1). Blowing up in Figure~\ref{cusp} yields Figure~\ref{5.2}. The handle slide drawn in Figure~\ref{ex-handleslide} gives Figure~\ref{5.3}. An additional blow-up yields Figure~\ref{5.4}, and an isotopy gives Figure~\ref{5.1}.

\begin{figure}[htbp]
\begin{center}
\begin{minipage}{.45\linewidth}
\begin{center}
\includegraphics[scale=0.75]{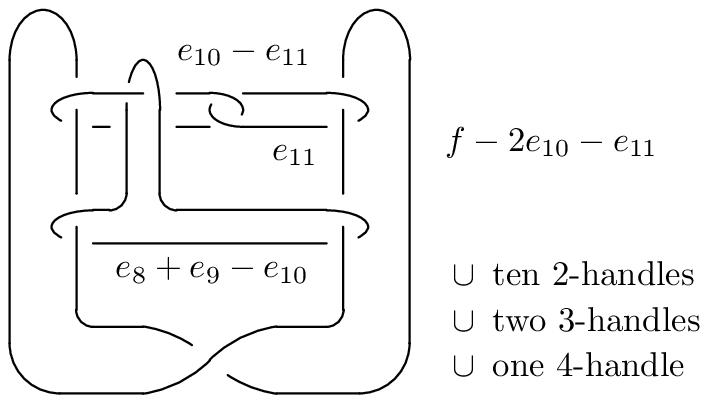}
\caption{$\mathbf{CP}^2\# 11\overline{\mathbf{C}\mathbf{P}^2}$}
\label{6.2}
\end{center}
\end{minipage}
\begin{minipage}{.45\linewidth}
\begin{center}
\includegraphics[scale=0.75]{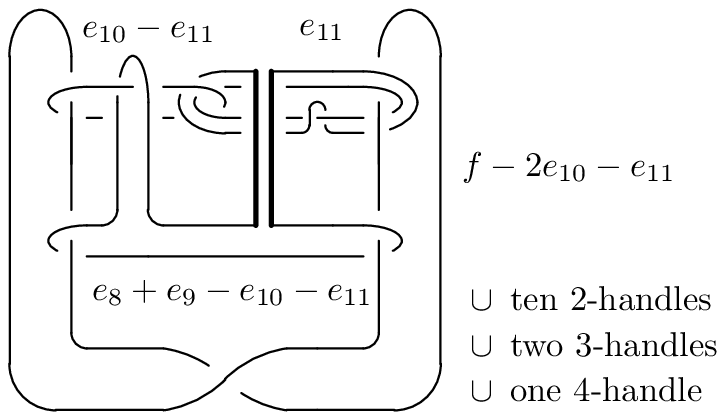}
\caption{$\mathbf{CP}^2\# 11\overline{\mathbf{C}\mathbf{P}^2}$}
\label{6.3}
\end{center}
\end{minipage}
\end{center}
\end{figure}

\begin{figure}[htbp]
\begin{center}
\includegraphics[scale=0.75]{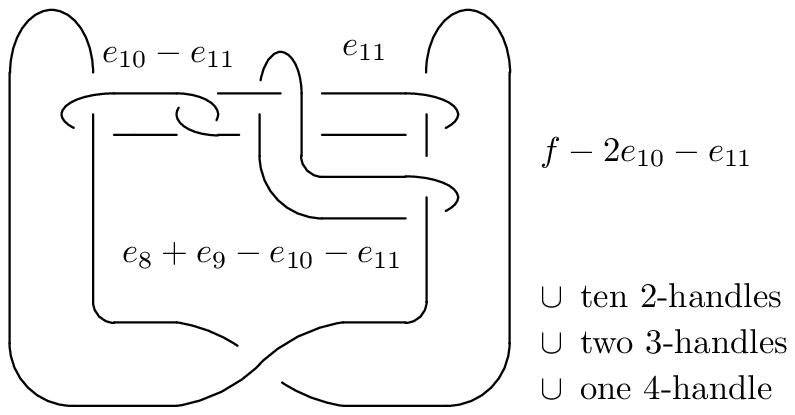}
\caption{$\mathbf{CP}^2\# 11\overline{\mathbf{C}\mathbf{P}^2}$}
\label{6.4}
\end{center}
\end{figure}

\begin{figure}[htbp]
\begin{center}
\includegraphics[scale=0.70]{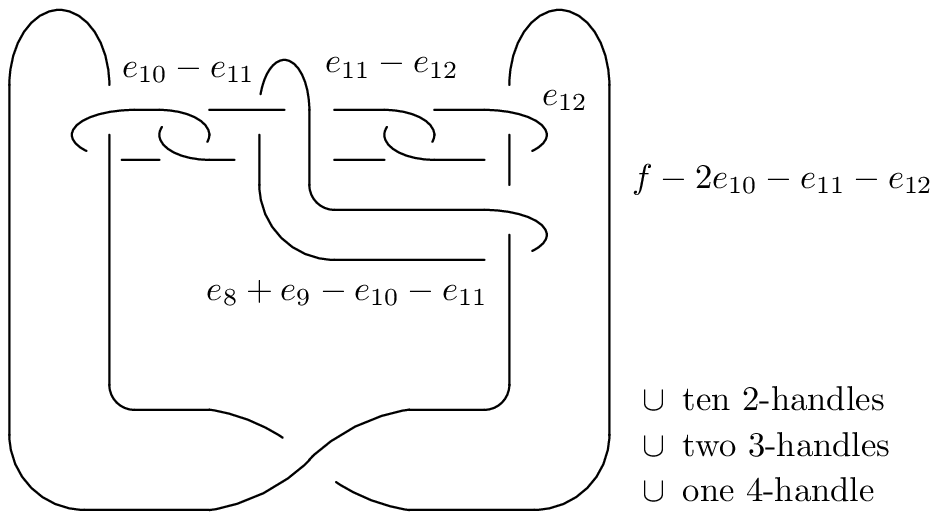}
\caption{$\mathbf{CP}^2\# 12\overline{\mathbf{C}\mathbf{P}^2}$}
\label{6.5}
\end{center}
\end{figure}

\begin{figure}[htbp]
\begin{center}
\includegraphics[scale=0.67]{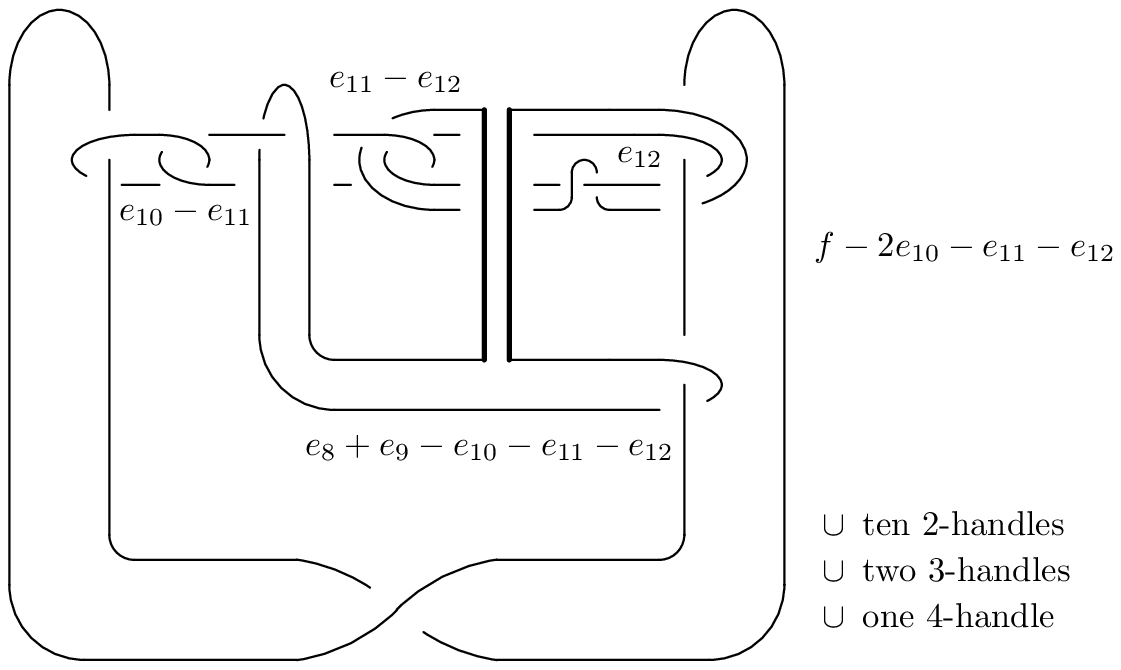}
\caption{$\mathbf{CP}^2\# 12\overline{\mathbf{C}\mathbf{P}^2}$}
\label{6.6.0}
\end{center}
\end{figure}

\begin{figure}[htbp]
\begin{center}
\includegraphics[scale=0.75]{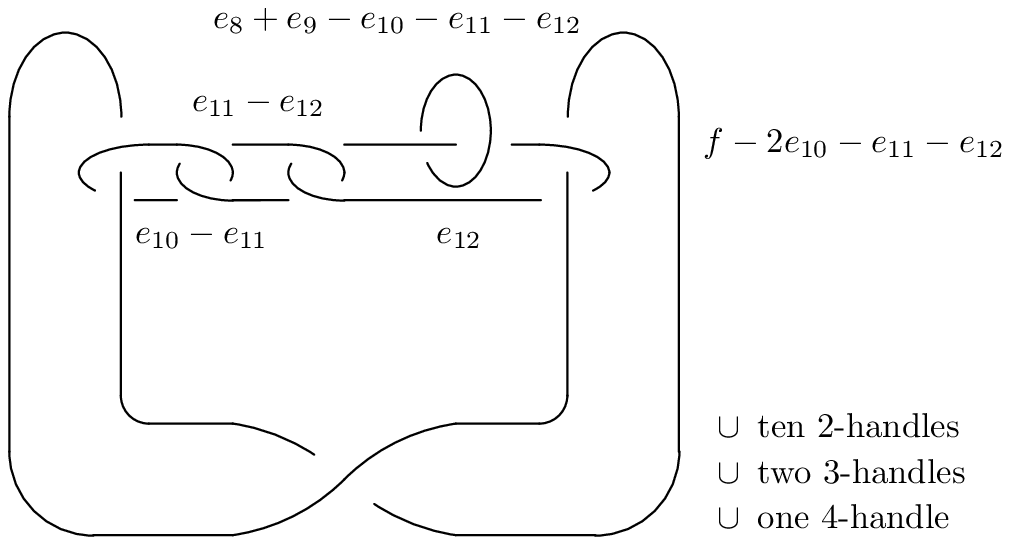}
\caption{$\mathbf{CP}^2\# 12\overline{\mathbf{C}\mathbf{P}^2}$}
\label{6.6}
\end{center}
\end{figure}

\begin{figure}[htbp]
\begin{center}
\includegraphics[scale=0.70]{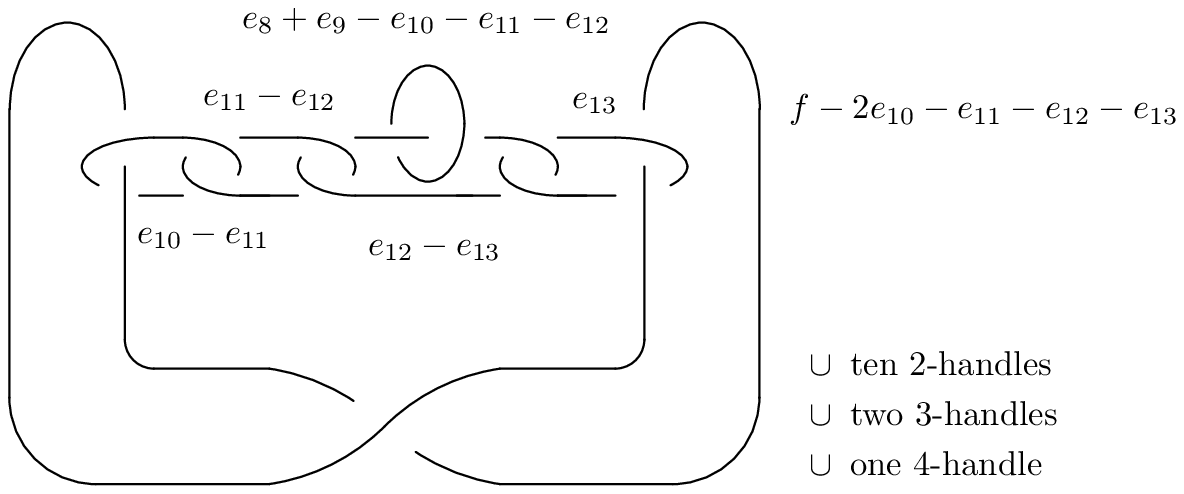}
\caption{$\mathbf{CP}^2\# 13\overline{\mathbf{C}\mathbf{P}^2}$}
\label{6.7}
\end{center}
\end{figure}

\begin{figure}[htbp]
\begin{center}
\begin{minipage}{.45\linewidth}
\begin{center}
\includegraphics[scale=0.80]{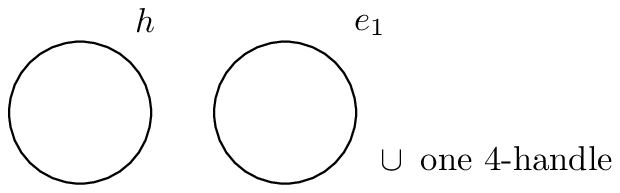}
\caption{$\mathbf{CP}^2\# \overline{\mathbf{C}\mathbf{P}^2}$}
\label{8.2}
\end{center}
\end{minipage}
\begin{minipage}{.45\linewidth}
\begin{center}
\includegraphics[scale=0.80]{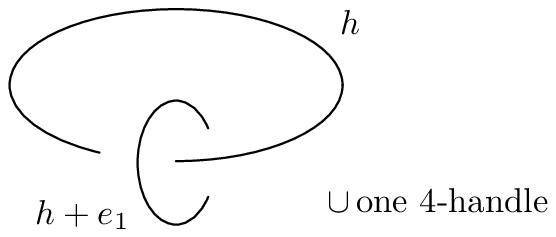}
\caption{$\mathbf{CP}^2\# \overline{\mathbf{C}\mathbf{P}^2}$}
\label{8.3}
\end{center}
\end{minipage}
\end{center}
\end{figure}

\begin{figure}[htbp]
\begin{center}
\begin{minipage}{.45\linewidth}
\begin{center}
\includegraphics[scale=0.80]{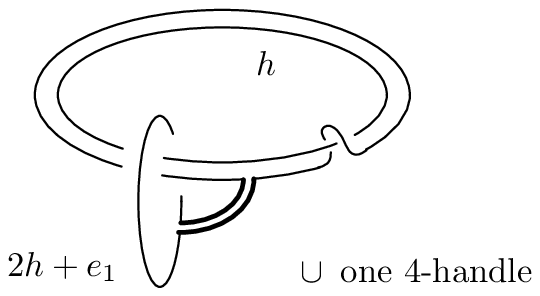}
\caption{$\mathbf{CP}^2\# \overline{\mathbf{C}\mathbf{P}^2}$}
\label{8.4}
\end{center}
\end{minipage}
\begin{minipage}{.45\linewidth}
\begin{center}
\includegraphics[scale=0.80]{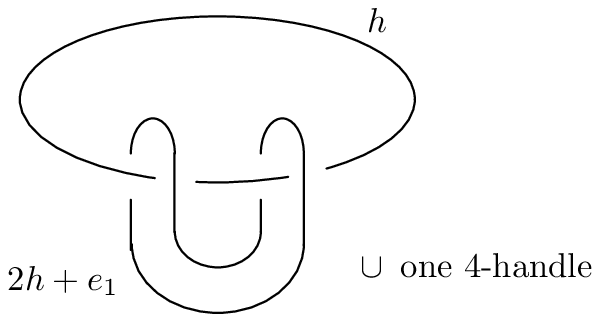}
\caption{$\mathbf{CP}^2\# \overline{\mathbf{C}\mathbf{P}^2}$}
\label{8.5}
\end{center}
\end{minipage}
\end{center}\medskip 

\begin{center}
\begin{minipage}{.45\linewidth}
\begin{center}
\includegraphics[scale=0.77]{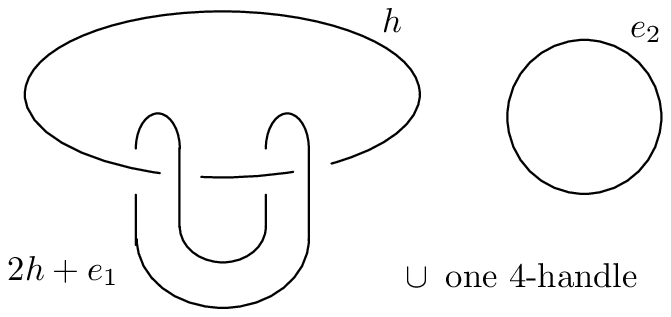}
\caption{$\mathbf{CP}^2\# 2\overline{\mathbf{C}\mathbf{P}^2}$}
\label{8.6}
\end{center}
\end{minipage}
\begin{minipage}{.45\linewidth}
\begin{center}
\includegraphics[scale=0.77]{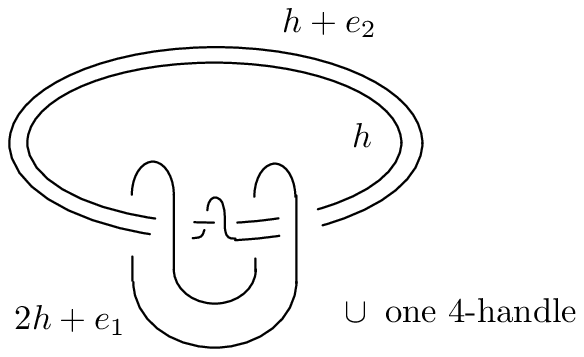}
\caption{$\mathbf{CP}^2\# 2\overline{\mathbf{C}\mathbf{P}^2}$}
\label{8.7}
\end{center}
\end{minipage}
\end{center}\medskip 

\begin{center}
\begin{minipage}{.45\linewidth}
\begin{center}
\includegraphics[scale=0.75]{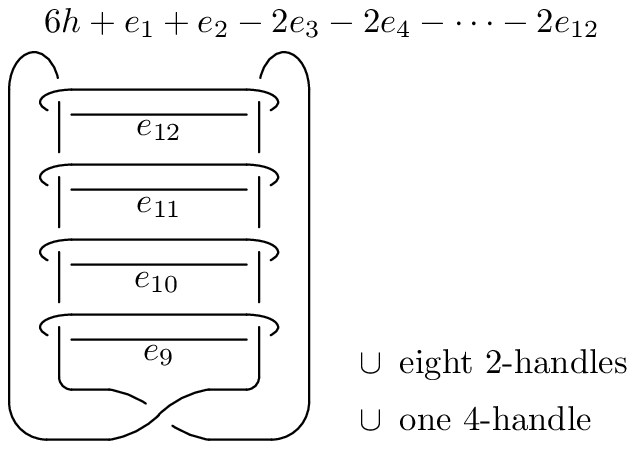}
\caption{$\mathbf{CP}^2\# 12\overline{\mathbf{C}\mathbf{P}^2}$}
\label{8.8}
\end{center}
\end{minipage}
\begin{minipage}{.45\linewidth}
\begin{center}
\includegraphics[scale=0.75]{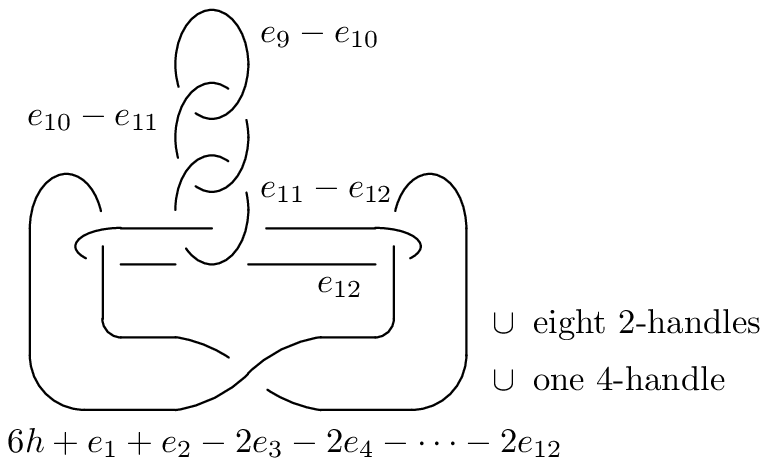}
\caption{$\mathbf{CP}^2\# 12\overline{\mathbf{C}\mathbf{P}^2}$}
\label{8.9}
\end{center}
\end{minipage}
\end{center}\medskip 

\begin{center}
\begin{minipage}{.45\linewidth}
\begin{center}
\includegraphics[scale=0.75]{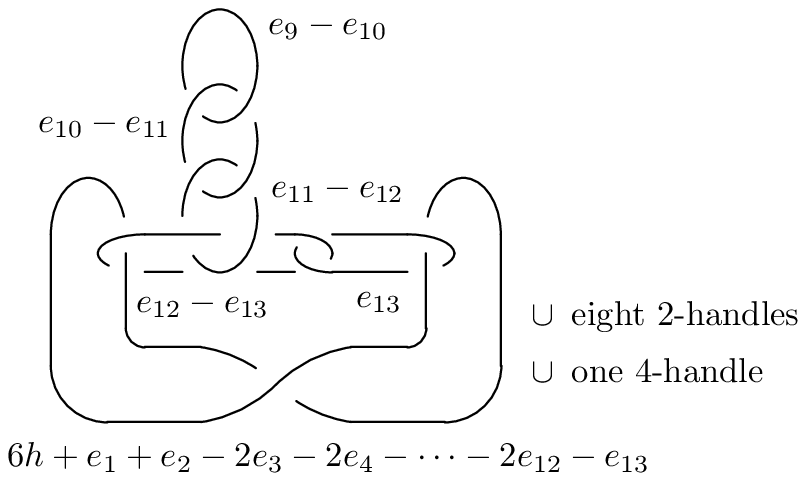}
\caption{$\mathbf{CP}^2\# 13\overline{\mathbf{C}\mathbf{P}^2}$}
\label{8.10}
\end{center}
\end{minipage}
\end{center}
\end{figure}

Secondly we give a proof for (2). Handle slides, isotopies and blow-ups in Figure~\ref{5.4} yield Figure~\ref{6.1}: 
\begin{multline*}
\text{Figure~\ref{5.4}} \xrightarrow{e_8+(e_9-e_{10})} \text{Figure~\ref{6.2}} 
\xrightarrow{\text{$(e_8+e_9-e_{10})-e_{11}$}} \text{Figure~\ref{6.3}}\\ 
\xrightarrow{\text{isotopy}} \text{Figure~\ref{6.4}} 
\xrightarrow{\text{blow-up}} \text{Figure~\ref{6.5}} 
\xrightarrow{\text{$(e_8+e_9-e_{10}-e_{11})-e_{12}$}} \text{Figure~\ref{6.6.0}}\\ 
\xrightarrow{\text{isotopy}} \text{Figure~\ref{6.6}} 
\xrightarrow{\text{blow-up}} \text{Figure~\ref{6.7}} 
\xrightarrow{\text{isotopy}} \text{Figure~\ref{6.1}}.\quad \rlap{\qedsymbol}
\end{multline*}
\begin{proposition}\label{prop:3.3}
$\mathbf{CP}^2\# 13\overline{\mathbf{C}\mathbf{P}^2}$ 
admits the handle decomposition drawn in Figure~\ref{8.1}. In particular $\mathbf{CP}^2\# 13\overline{\mathbf{C}\mathbf{P}^2}$ contains the copy of $C_5$ drawn in the figure.
\begin{figure}[ht!]
\begin{center}
\includegraphics[scale=0.75]{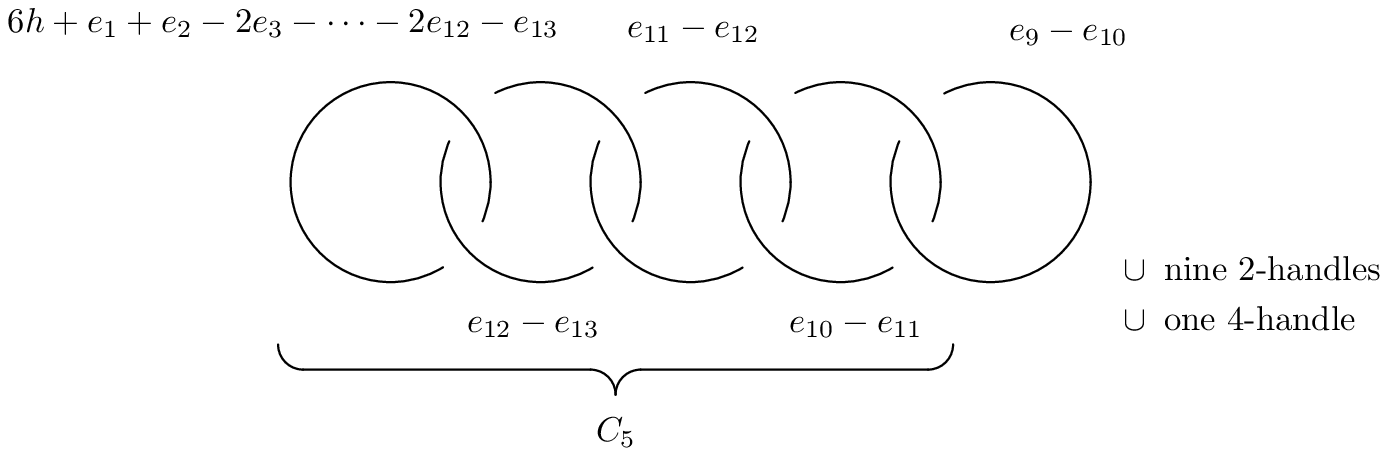}
\caption{$\mathbf{CP}^2\# 13\overline{\mathbf{C}\mathbf{P}^2}$}
\label{8.1}
\end{center}
\end{figure}
\end{proposition}
\begin{proof}Recall the construction in the proof of Lemma~\ref{lem:3.1}. In this construction, we created a $2$-handle/$3$-handle pair twice. Instead of introducing a $2$-handle/$3$-handle pair twice, blowing up twice yields Figure~\ref{8.7}:
\begin{multline*}
\quad \text{Figure~\ref{4.1}} \xrightarrow{\text{blow-up}} \text{Figure~\ref{8.2}} \xrightarrow{\text{$e_1+h$}} \text{Figure~\ref{8.3}} \xrightarrow{\text{$(h+e_1)+h$}} \text{Figure~\ref{8.4}}\\ \xrightarrow{\text{isotopy}} \text{Figure~\ref{8.5}} \xrightarrow{\text{blow-up}} \text{Figure~\ref{8.6}} \xrightarrow{\text{$e_2+h$}} \text{Figure~\ref{8.7}}.\quad 
\end{multline*}
Handle slides and blow-ups as in proofs of Lemma~\ref{lem:3.1} and Proposition~\ref{prop:3.2} gives Figure~\ref{8.8}. Repeating handle slides drawn in Figure~\ref{ex-handleslide} yields Figure~\ref{8.9}. An additional blow-up gives Figure~\ref{8.10}, and an isotopy gives Figure~\ref{8.1}.
\end{proof}
\begin{definition}
Let $E_q$ be the smooth $4$-manifold obtained from $\mathbf{CP}^2\# (8+q)\overline{\mathbf{C}\mathbf{P}^2}$ by rationally blowing down along the copy of $C_q$ obtained in Proposition~\ref{prop:3.2}, for $q=3,5$. 
Let $E'_3$ be the smooth $4$-manifold obtained from $\mathbf{CP}^2\# 13\overline{\mathbf{C}\mathbf{P}^2}$ by rationally blowing down along the copy of $C_5$ obtained in Proposition~\ref{prop:3.3}.
\end{definition}
\begin{remark}\normalfont 
It is not known whether or not there exists a copy of $C_5$ in $\mathbf{CP}^2\# 13\overline{\mathbf{C}\mathbf{P}^2}$ such that the rational blow-down is diffeomorphic to $E(1)_{2,3}$.

In \cite{Y3} we will construct more examples of exotic $\mathbf{CP}^2\# 9\overline{\mathbf{C}\mathbf{P}^2}$ without $1$- and $3$-handles, by improving the construction of $E_3'$. The author does not know if these examples have the same Seiberg-Witten invariants as the elliptic surfaces $E(1)_{p,q}$. 
\end{remark}
We prepare the following lemma.
\begin{lemma}[{cf.~Gompf-Stipsicz \cite{GS}}]\label{without-handle}Suppose that a simply connected closed smooth 
$4$-manifold $X$ has the handle decomposition drawn 
in Figure~\ref{without}. Here $n$ is an arbitrary integer, $h_2$ and $h_3$ are arbitrary natural numbers. Note that we write usual framings instead of homology classes in the figure. 

Let $X_{(p)}$ be the rational blow-down of $X$ along the copy of $C_p$ drawn in Figure~\ref{without}. 
Then $X_{(p)}$ admits a handle decomposition 
\begin{equation*}
X_{(p)}=\text{one $0$-handle} \cup \text{$(h_2+1)$ $2$-handles} \cup \text{$h_3$ $3$-handles} \cup \text{one $4$-handle}.
\end{equation*}
In particular $X_{(p)}$ admits a handle decomposition without $1$-handles.
\begin{figure}[ht!]
\begin{center}
\includegraphics[scale=1.0]{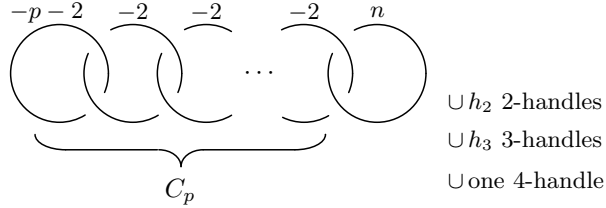}
\caption{Handle decomposition of $X$}
\label{without}
\end{center}
\end{figure}
\end{lemma}
\begin{proof}Draw a Kirby diagram of $X_{(p)}$, following the procedure introduced in \cite[Section 8.5]{GS} (see also \cite[page 516 Solution of Exercise 8.5.1.(a)]{GS}). Then the $n$-framed unknot drawn in Figure~\ref{without} changes into a meridian of a unique dotted circle which naturally appears in this procedure. Thus we can cancel the $1$-handle/$2$-handle pair. Note that this procedure does not produce new $3$-handles.
\end{proof}
The following proposition gives Theorem 1.1.(1)(a)(c) and (2)(a)(c).
\begin{proposition}
For $q=3,5$, the manifold $E_q$ is homeomorphic to $E(1)_{2,q}$ and admits a handle decomposition without $1$-handles, namely,\begin{equation*}
E_q=\text{one $0$-handle} \cup \text{twelve $2$-handles} \cup \text{two $3$-handles} \cup \text{one $4$-handle}.
\end{equation*}

$E_3'$ is homeomorphic to $E(1)_{2,3}$ and admits a handle decomposition without $1$- and $3$-handles, namely,
\begin{equation*}
E_3'=\text{one $0$-handle} \cup \text{ten $2$-handles} \cup \text{one $4$-handle}.
\end{equation*}
\end{proposition}
\begin{proof}
Lemma~\ref{without-handle} shows the above properties of $E_q$ and $E_3'$ about handle decompositions. Thus $E_q$ and $E_3'$ are simply connected. Since $E_q$ is obtained from $\mathbf{CP}^2\# (8+q)\overline{\mathbf{C}\mathbf{P}^2}$ by rationally blowing down along a copy of $C_q$, we have 
\begin{align*}
 b_2^+(E_q)&=b_2^+(\mathbf{CP}^2\# (8+q)\overline{\mathbf{C}\mathbf{P}^2})=1,\\
 b_2^-(E_q)&=b_2^-(\mathbf{CP}^2\# (8+q)\overline{\mathbf{C}\mathbf{P}^2})-b_2^-(C_q)=(8+q)-(q-1)=9.\end{align*}
Similarly we have $b_2^+(E_3')=1$ and $b_2^-(E_3')=9$. Therefore Freedman's theorem together with Rochlin's theorem shows that $E_q$ and $E_3'$ are homeomorphic to $\mathbf{CP}^2\# 9\overline{\mathbf{C}\mathbf{P}^2}$. Thus $E_q$ and $E_3'$ are homeomorphic to $E(1)_{2,q}$. 
\end{proof}
\section{Seiberg-Witten invariants}
In this section, we briefly review facts about the Seiberg-Witten invariants with $b_2^+=1$. 
For details and examples of computations, see Fintushel-Stern \cite{FS3}, \cite{FS1}, \cite{FS2}, Stern \cite{S}, Park \cite{P1}, \cite{P2}, Ozsv\'{a}th-Szab\'{o} \cite{OS}, Stipsicz-Szab\'{o} \cite{SS} and Park-Stipsicz-Szab\'{o} \cite{PSS}.

Suppose that $X$ is a simply connected closed smooth $4$-manifold with $b_2^+(X)=1$. Let $\mathcal{C}(X)$ be the set of 
characteristic elements of $H^2(X;\mathbf{Z})$. Fix a homology orientation on $X$, that is, 
orient $H^2_+(X;\mathbf{R}):=\{ H\in H^2(X;\mathbf{Z})\, |\, H^2>0\} $. 
Then the (small-perturbation) Seiberg-Witten invariant $SW_{X,H}(K)\in \mathbf{Z}$ is defined for 
every positively oriented element $H\in H^2_+(X;\mathbf{R})$ and every element $K\in \mathcal{C}(X)$ such that $K\cdot H\neq 0$. 
Let $e(X)$ and $\sigma(X)$ be the Euler characteristic and the signature of $X$, respectively, and $d_X(K)$ the even integer defined by 
$d_X(K)=\frac{1}{4}(K^2-2e(X)-3\sigma(X))$ 
for $K\in \mathcal{C}(X)$. It is known that if $SW_{X,H}(K)\neq 0$ for some $H\in H^2_+(X;\mathbf{R})$, then $d_X(K)\ge 0$. 
The wall-crossing formula tells us the dependence of $SW_{X,H}(K)$ on $H$:
if $H, H' \in H^2_+(X;\mathbf{R})$ and $K\in \mathcal{C}(X)$ satisfy $H\cdot H'>0$ and $d_X(K)\ge 0$, then
\begin{multline*}
SW_{X,H'}(K)=SW_{X,H}(K)\\
+
\begin{cases}
0&\text{if $K\cdot H$ and $K\cdot H'$ have the same sign,}\\
(-1)^{\frac{1}{2}d_X(K)}&\text{if $K\cdot H>0$ and $K\cdot H'<0$,}\\
(-1)^{1+\frac{1}{2}d_X(K)}&\text{if $K\cdot H<0$ and $K\cdot H'>0$}.
\end{cases}
\end{multline*}
Note that these facts imply that $SW_{X,H}(K)$ is independent of $H$ in the case $b_2^-(X)\leq 9$, in other words, 
the Seiberg-Witten invariant $SW_{X}:\mathcal{C}(X)\to \mathbf{Z}$ is well-defined.

We recall the change of the Seiberg-Witten invariants by rationally blowing down. Assume that $X$ contains a copy of $C_p$. Let $X_{(p)}$ be the rational blow-down of $X$ along the copy of $C_p$. Suppose that $X_{(p)}$ is simply connected. The following theorems are known.
\begin{proposition}[{Fintushel-Stern \cite{FS1}}]\label{thm:4.1}
For every element $K\in \mathcal{C}(X_{(p)})$, there exists an element 
$\tilde{K}\in \mathcal{C}(X)$ such that 
$K\rvert _{X_{(p)}-\text{\normalfont{int}}\,B_{p}}=\tilde{K}\rvert _{X-\text{\normalfont{int}}\,C_{p}}$ and 
$d_{X_{(p)}}(K)=d_X(\tilde{K})$. We call such an element $\tilde{K}\in \mathcal{C}(X)$ a {\normalfont \textrm{lift}} of $K$.
\end{proposition}
\begin{theorem}[{Fintushel-Stern \cite{FS1}}]\label{nthm:4.1}
If an element $\tilde{K}\in \mathcal{C}(X)$ is a lift of some element $K\in \mathcal{C}(X_{(p)})$, then $SW_{X_{(p)},H}(K)=SW_{X,H}(\tilde{K})$ for every positively oriented element $H\in H^2_+(X;\mathbf{R})$ 
which is orthogonal to the subspace $H_2(C_p;\mathbf{R})$ of $H_2(X;\mathbf{R})$. Note that we view $H$ as a positively oriented element of $H^2_+(X_{(p)};\mathbf{R})$.
\end{theorem}
\begin{theorem}[{Fintushel-Stern \cite{FS1}, cf.~Park \cite{P1}}]\label{thm:4.2}
If an element $\tilde{K}\in \mathcal{C}(X)$ satisfies that $(\tilde{K}\rvert _{C_p})^2=1-p$ and 
$\tilde{K}\rvert _{\partial C_p}=mp\in \mathbf{Z}_{p^2}\cong H^2(\partial C_p;\mathbf{Z})$ 
with $m\equiv p-1\pmod 2$, then there exists an element 
$K\in \mathcal{C}(X_{(p)})$ such that $\tilde{K}$ is a lift of $K$. 
\end{theorem}
\begin{corollary}\label{cor:4.3}
If an element $\tilde{K}\in \mathcal{C}(X)$ satisfies $\tilde{K}(u_1)=\dots=\tilde{K}(u_{p-2})=0$ and $\tilde{K}(u_{p-1})=\pm p$, then $\tilde{K}$ is a lift of some element $K\in \mathcal{C}(X_{(p)})$.
\end{corollary}
\section{Computations of SW invariants}
In this section we complete the proof of Theorem 1.1. We prepare the following lemma here. 
\begin{lemma}\label{lem:5.1}
Let $X$ be a simply connected closed smooth $4$-manifold which contains a copy of $C_p$, and $\iota$ the inclusion $X-\text{\normalfont{int}}\,C_p\hookrightarrow X$. Let $C_p^{\perp }$ be the orthogonal complement of the subspace spanned by $u_1,\dots,u_{p-1}\in H_2(X;\mathbf{Z})$, that is, 
\begin{equation*}
C_p^{\perp }:=\{ v\in H_2(X;\mathbf{Z})\,|\,v\cdot u_1=\dots=v\cdot u_{q-1}=0\}.
\end{equation*}

Suppose that there exists an element $\delta \in H_2(X;\mathbf{Z})$ such that $\delta \cdot u_1=1$ and $\delta \cdot u_2=\delta \cdot u_3=\dots=\delta \cdot u_{p-1}=0$. 
Then 
\begin{align}
 &\iota_*H_2(X-\text{\normalfont{int}}\,C_p;\mathbf{Z})=C_p^{\perp };\\
 &H_1(X-\text{{\normalfont int}}\,C_p;\mathbf{Z})=0.
\end{align}
\end{lemma}
\begin{proof}
Firstly we give a proof for (1). 
Since every element of $H_2(X-\text{int}\,C_p;\mathbf{Z})$ is represented by a surface, it is clear that $\iota_*H_2(X-\text{int}\,C_p;\mathbf{Z})\subset C_p^{\perp }$. 

Let $\iota '$ be the inclusion $C_p\hookrightarrow X$. Mayer-Vietoris exact sequence of $(X-\text{int}\,C_p)\cup C_p=X$ is as follows:
\begin{equation*}
0\to H_2(X-\text{int}\,C_p;\mathbf{Z})\oplus H_2(C_p;\mathbf{Z})\stackrel{\iota_*+ \iota '_*}{\to }
H_2(X;\mathbf{Z})\stackrel{\partial}{\to } \mathbf{Z}_{p^2}.
\end{equation*}
Since $C_p$ is negative definite and $\text{Im}\,\iota _*\subset C_p^{\perp }$, we have $\text{Im}\,(\iota _*+\iota '_*)=\text{Im}\,\iota _*\oplus \text{Im}\,\iota '_*$. 

We determine $\partial (\delta)$ here. 
There clearly exists an element $n\in \mathbf{Z}$ such that $\partial (n\delta )\equiv 0\pmod{p^2}$. The above exact sequence ensures the existence of elements $u\in \text{Im}\,\iota '_*$ and $v\in C_p^\bot $ such that $n\delta =u+v$. 
The element $u$ satisfies $u\cdot u_1=n$ $(=n\delta \cdot u_1)$ and $u\cdot u_2=u\cdot u_3=\dots=u\cdot u_{p-1}=0$ $(=n\delta \cdot u_2)$. 
Since $u_1,u_2,\dots,u_n$ is a basis of $\text{Im}\,\iota '_*$, we can easily see $n \equiv 0\pmod{p^2}$ by using the intersection form of $C_p$. 
Hence $\partial (\delta)$ is a generator of $\mathbf{Z}_{p^2}$.

Suppose that some element $w\in C_p^\bot $ satisfies $\partial(w)\not\equiv 0\pmod{p^2}$. 
Since $\partial (\delta)$ is a generator of $\mathbf{Z}_{p^2}$, there exists an element $n'\in \mathbf{Z}$ with $n' \not\equiv 0\pmod{p^2}$ such that $\partial (n'\delta +w)\equiv 0$. 
Applying the above argument about $n\delta$ to $n'\delta +w$, we get $n'\equiv 0\pmod{p^2}$. This is a contradiction. Thus we obtain $\partial (C_p^\bot )=0$. Therefore $C_p^\bot \subset \text{Ker}\,\partial =\text{Im}\,\iota _*\oplus \text{Im}\,\iota '_*\subset C_p^\bot \oplus \text{Im}\,\iota '_*$. Thus it is easy to see $C_p^\bot \subset \iota_*H_2(X-\text{int}\,C_p;\mathbf{Z})$. 

Secondly we give a proof for (2). Since the above $\partial$ is onto, we can easily show by using Mayer-Vietoris exact sequence. 
\end{proof}
\begin{remark}
(1) Since $\iota_*: H_2(X-\text{\normalfont{int}}\,C_p;\mathbf{Z})\to H_2(X;\mathbf{Z})$ is injective, the above lemma allows us to identify $H_2(X-\text{\normalfont{int}}\,C_p;\mathbf{Z})$ with $C_p^{\perp }$. \\
(2) Under the same assumption as that in Lemma~\ref{lem:5.1}, we can also show $H_1(X_{(p)};\mathbf{Z})=0$. Here $X_{(p)}$ denotes the rational blow-down of $X$ along the copy of $C_p$. It is not known whether or not the fundamental groups of $X-\text{{\normalfont int}}\,C_p$ and $X_{(p)}$ vanish. 
\end{remark}

The following proposition gives us Theorem 1.1.(1)(b). In the rest of this section, we denote the symbol $R_n$ as $\mathbf{CP}^2\# n\overline{\mathbf{C}\mathbf{P}^2}$. 
\begin{proposition}\label{prop:5.2}
$E_q$ has the same Seiberg-Witten invariant as $E(1)_{2,q}$, that is, there exists a homeomorphism between $E_q$ and $E(1)_{2,q}$ which preserves the orientations, the homology orientations and the Seiberg-Witten invariants, for $q=3,5$.
\end{proposition}
\begin{proof}
We give a proof for $q=3$, firstly. 
Let $\alpha_1,\alpha_2,\dots,\alpha_9,\beta\in {_2}C_3^{\perp}$ be the elements defined by
\begin{align*}
\alpha_1 &= 4h-e_1-e_2-\dots-e_9-2e_{10}-2e_{11},\\
\alpha_i &= 5h-2e_1-2e_2-e_3-e_4-\dots-e_9-2e_{10}-2e_{11}-e_{i+1}\, (2\leq i\leq 8),\\
\alpha_9 &= e_1-e_2,\: 
\beta = 30h-13e_1-10e_2-7e_3-7e_4-\dots-7e_9-12e_{10}-12e_{11}.
\end{align*}
We can view $\alpha_1,\alpha_2,\dots,\alpha_9,\beta$ as elements of $H_2(E(1)_{2,3};\mathbf{Z})$ by Lemma~\ref{lem:5.1}.(1), Corollary~\ref{cor:2.2} 
and the following natural identification:
\begin{equation*}
H_2(E(1)_{2,3}-\text{int}\,B_3;\mathbf{Z})
(=H_2(R_{11}-\text{int}\,C_3;\mathbf{Z}))
\subset H_2(E(1)_{2,3};\mathbf{Z}).
\end{equation*}
This identification preserves cup products. Therefore the elements $\alpha_1,\alpha_2,\dots,
\alpha_9,\beta $ of $H_2(E(1)_{2,3};\mathbf{Z})$ satisfy 
\begin{align*}
 &\alpha_1^2=\alpha_2^2=\dots=\alpha_8^2=-1,\, \alpha_9^2=-2,\, \alpha_i\cdot \alpha_j=0\, (1\leq i<j\leq 9),\\
&\beta^2=0,\, \beta\cdot \alpha_1=\beta\cdot \alpha_2=\dots=\beta\cdot \alpha_8=0,\, \beta\cdot \alpha_9=3.
\end{align*}

Recall that the intersection form of $E(1)_{2,3}$ is $\langle 1 \rangle \oplus 9\langle -1 \rangle$ (This notation of the intersection form is the same as that in Gompf-Stipsicz \cite[Section 1.2]{GS}.). This implies that either the matrix 
$\left(
\begin{smallmatrix}
0 & 1\\
1 & 0
\end{smallmatrix}
\right)$
 or 
$\left(
\begin{smallmatrix}
1 & 0\\
0 & -1
\end{smallmatrix}
\right)$
represents the symmetric bilinear form on $\langle \alpha_1,\alpha_2,\dots,\alpha_8 \rangle^{\bot }$. We here denote the symbol $\langle \alpha_1,\alpha_2,\dots,\alpha_8 \rangle^{\bot }$ as the orthogonal complement of the subspace spanned by $\alpha_1,\alpha_2,\dots,\alpha_8\in H_2(E(1)_{2,3};\mathbf{Z})$. 
Since $\alpha_9$ and $\beta $ are elements of $\langle \alpha_1,\alpha_2,\dots,\alpha_8 \rangle^{\bot }$, it is easy to check that the matrix
$\left(
\begin{smallmatrix}
0 & 1\\
1 & 0
\end{smallmatrix}
\right)$ 
represents the symmetric bilinear form on $\langle \alpha_1,\alpha_2,\dots,\alpha_8 \rangle^{\bot }$. 
We can easily see that there exists an element $\alpha_{10}\in H_2(E(1)_{2,3};\mathbf{Z})$ such that $3\alpha_{10}=\beta$, by using a basis of $\langle \alpha_1,\alpha_2,\dots,\alpha_8 \rangle^{\bot }$. Note that $\alpha_1,\alpha_2,\dots,\alpha_{10}$ is a basis of $H_2(E(1)_{2,3};\mathbf{Z})$. 

Proposition 3.2.(1) allows us to apply the above argument to $E_3$. Thus we get a basis $\alpha_1',\dots,\alpha_{10}'$ of $H_2(E_3;\mathbf{Z})$ which is corresponding to the basis $\alpha_1,\dots,\alpha_{10}$ of $H_2(E(1)_{2,3};\mathbf{Z})$. Let $\varphi:H^2(E(1)_{2,3};\mathbf{Z})\to H^2(E_3;\mathbf{Z})$ be the isomorphism defined by $PD(\alpha_i)\mapsto PD(\alpha_i')$ $(1\leq i\leq 10)$. Here $PD$ denotes Poincar\'{e} dual. The isomorphism $\varphi$ preserves the intersection forms and the homology orientations of $E(1)_{2,3}$ and $E_3$.

Proposition~\ref{thm:4.1} gives us a lift $\tilde{K}\in \mathcal{C}(R_{11})$ of $K$ for every $K\in \mathcal{C}(E(1)_{2,3})$. Lemma~\ref{lem:5.1}.(2) together with the universal coefficient theorem implies that $\tilde{K}\rvert _{R_{11}-\text{\normalfont{int}}\,C_{3}}$ and $\varphi(K)\rvert _{E_3-\text{\normalfont{int}}\,B_{3}}$ are uniquely determined by their values on $H_2(R_{11}-\text{\normalfont{int}}\,C_{3};\mathbf{Z})=H_2(E_3-\text{\normalfont{int}}\,B_{3};\mathbf{Z})$. Since $\alpha_i'=\alpha_i$ $(1\leq i\leq 9)$ and $3\alpha_{10}'=3\alpha_{10}$ as elements of $_2C_3^{\perp}$, it is easy to check that $\tilde{K}$ is also a lift of the element $\varphi(K)\in \mathcal{C}(E_3)$. Thus Theorem~\ref{nthm:4.1} shows $SW_{E_3}(\varphi(K))=SW_{E(1)_{2,3}}(K)$. Hence the isomorphism $\varphi$ preserves the Seiberg-Witten invariants of $E(1)_{2,3}$ and $E_3$. 
Freedman's theorem gives us a required homeomorphism $\Phi:E_3\to E(1)_{2,3}$ which preserves the orientations and satisfies $\Phi^*=\varphi$.

We briefly give a proof for $q=5$, secondly. 
Let $\alpha_{5,1},\alpha_{5,2},\dots,\alpha_{5,9},\beta_5\in {_2}C_5^{\perp}$ be the elements defined by
\begin{align*}
\alpha_{5,i} &= 17h-3e_1-4e_2-\dots-4e_9-6e_{10}-\dots-6e_{13}-e_{i+1}\: (1\leq i\leq 8),\\
\alpha_{5,9} &= 96h-19e_1-23e_2-\dots-23e_9-34e_{10}-\dots-34e_{13},\\
\beta_5 &= 537h-104e_1-129e_2-\dots-129e_9-190e_{10}-\dots-190e_{13}.
\end{align*}
Applying the above argument to $E_5$, we obtain a proof.
\end{proof}

To prove Theorem 1.1.(2)(b), we compute the Seiberg-Witten invariant of $E(1)_{2,3}$. 
\begin{lemma}\label{lem:5.3}
Let $K_3\in \mathcal{C}(E(1)_{2,3})$ be the element defined by $K_3=PD(\alpha_1+\dots+\alpha_8-2\alpha_9-4\alpha_{10})$. Here $\alpha _1,\alpha _2,\dots,\alpha _{10}$ denote the elements of $H_2(E(1)_{2,3};\mathbf{Z})$ defined in the proof of Proposition~\ref{prop:5.2}. Then $K_3$ satisfies $SW_{E(1)_{2,3}}(\pm K_3)=\pm 1$ and is the unique element of $\mathcal{C}(E(1)_{2,3})$ up to sign for which $SW_{E(1)_{2,3}}$ is nonzero.
\end{lemma}
\begin{proof}
Let $\tilde{K_3}\in \mathcal{C}(R_{11})$ and $H\in H_+^2(R_{11};\mathbf{R})$ be the elements defined by $\tilde{K_3}=PD(3h-e_1-e_2-\dots-e_{11})$ and $H=PD(7h-2e_1-2e_2-\dots-2e_{11})$. Note that $H$ is orthogonal to the subspace $H_2(C_3;\mathbf{R})$ of $H_2(R_{11};\mathbf{R})$. It is well known that $SW_{R_{n},PD(h)}(\tilde{K})=0$ for every $\tilde{K}\in \mathcal{C}(R_{n})$ and every $n\geq 0$.
Applying the wall-crossing formula to $\pm \tilde{K_3},\, H$ and $PD(h)$, we have $SW_{R_{11},H}(\pm \tilde{K_3})=\pm 1$. 
Corollary~\ref{cor:4.3} shows that $\tilde{K_3}$ is a lift of some element $K_3\in \mathcal{C}(E(1)_{2,3})$. Thus Theorem~\ref{nthm:4.1} gives $SW_{E(1)_{2,3}}(\pm K_3)=\pm 1$. 

Since $\tilde{K_3}$ is a lift of $K_3$, the element $K_3$ satisfies $K_3(\alpha_i)=\tilde{K}_3(\alpha_i)$ $(1\leq i\leq 9)$ and $K_3(3\alpha_{10})=\tilde{K}_3(3\alpha_{10})$. 
Hence the values of $K_3$ are as follows: $K_3(\alpha_1)=K_3(\alpha_2)=\dots=K_3(\alpha_8)=-1$, $K_3(\alpha_9)=0$ and $K_3(\alpha_{10})=-2$. Therefore we get $K_3=PD(\alpha_1+\dots+\alpha_8-2\alpha_9-4\alpha_{10})$.

Suppose that an element $L\in \mathcal{C}(E(1)_{2,3})$ satisfies $SW_{E(1)_{2,3}}(L)\neq 0$. 
Proposition~\ref{thm:4.1} ensures the existence of a lift $\tilde{L}\in \mathcal{C}(R_{11})$ of $L$ such that $SW_{R_{11},H}(\tilde{L})\neq 0$. 
We put $a:=\tilde{L}(h)$. Since $L$ is characteristic and $d_{R_{11}}(\tilde{L})\geq 0$, the integer $a$ is odd and $\lvert a\rvert \geq 3$. In the case $a\geq 3$, Cauchy-Schwartz inequality ($(x_1y_1+\dots+x_ny_n)^2\leq (x_1^2+\dots+x_n^2)(y_1^2+\dots+y_n^2)$ for $x_1,\dots,x_n,y_1,\dots,y_n\in \mathbf{R}$) and $d_{R_{11}}(\tilde{L})=\frac{1}{4}(a^2-((\tilde{L}(e_1))^2+(\tilde{L}(e_2))^2+\dots+(\tilde{L}(e_{11}))^2)+2)\geq 0$ show 
\begin{align*}
\tilde{L}\cdot H&=7a-2\tilde{L}(e_1)-2\tilde{L}(e_2)-\dots-2\tilde{L}(e_{11})\\
&\geq 7a-\sqrt{2^2+2^2+\dots+2^2}\sqrt{(\tilde{L}(e_1))^2+(\tilde{L}(e_2))^2+\dots+(\tilde{L}(e_{11}))^2}\\
 &\geq 7a-2\sqrt{11}\sqrt{a^2+2}.
\end{align*}
Since $SW_{E(1)_{2,3}}(L)\neq 0$ and $a\geq 3$, the wall-crossing formula shows $\tilde{L}\cdot H<0$. Therefore we get $a=3$. This together with $\tilde{L}\cdot H<0$ shows $\tilde{L}(e_i)=1$ $(1\leq i\leq 11)$. We thus have $\tilde{L}=\tilde{K_3}$. 
Similarly we have $\tilde{L}=-\tilde{K_3}$ in the case $a\leq -3$. Hence $L=\pm K_3$.
\end{proof}
The following proposition completes the proof of Theorem~1.1.
\begin{proposition}\label{prop:5.4}
$E_3'$ has the same Seiberg-Witten invariant as $E(1)_{2,3}$, that is, there exists a homeomorphism between $E_3'$ and $E(1)_{2,3}$ which preserves the orientations, the homology orientations and the Seiberg-Witten invariants. 
\end{proposition}
\begin{proof}
Let $\tilde{K}_3'\in \mathcal{C}(R_{13})$ and $H'\in H_+^2(R_{13};\mathbf{R})$ be the elements defined by $\tilde{K}_3'=PD(3h+e_1+e_2-e_3-\dots-e_{13})$ and $H'=PD(23h+6e_1+6e_2-6e_3-\dots-6e_{13})$. Note that $H'$ is orthogonal to the subspace $H_2(C_5;\mathbf{R})$ of $H_2(R_{13};\mathbf{R})$. 
Applying the wall-crossing formula to $\pm {\tilde{K}_3}',\, H'$ and $PD(h)$, we get $SW_{R_{13},H'}(\pm \tilde{K}_3')=\pm 1$. Corollary~\ref{cor:4.3} shows that $\tilde{K}_3'$ is a lift of some element $K_3'\in \mathcal{C}(E_3')$. Thus Theorem~\ref{nthm:4.1} gives $SW_{E_3'}(\pm K_3')=\pm 1$. The same argument as that in the proof of Lemma~\ref{lem:5.3} shows that $K_3'$ is the unique element up to sign for which $SW_{E'_3}$ is nonzero.

Let $\alpha'\in H_2(R_{13};\mathbf{Z})$ be the element defined by $\alpha'=3h+e_1-e_3-e_4-\dots-e_7-e_{10}-e_{11}-e_{12}-e_{13}$. Lemma~\ref{lem:5.1}.(1) allows us to view $\alpha'$ as an element of $H_2(E_3';\mathbf{Z})$. We set $L_3'\in H^2(E_3';\mathbf{Z})$ by $L_3'=K_3'-PD(\alpha')$. The element $L_3'$ is a characteristic element of $\langle PD(\alpha') \rangle ^{\bot }$ and satisfies ${L_3'}^2=1$ and $K_3'=L_3'+PD(\alpha')$, because ${K_3'}^2=0$, $K_3\cdot PD(\alpha')=-1$ and $(PD(\alpha'))^2=-1$. We here denote the symbol $\langle PD(\alpha') \rangle ^{\bot }$ as the orthogonal complement of the subspace spanned by $PD(\alpha')\in H^2(E_3';\mathbf{Z})$. Since the symmetric bilinear form on $\langle PD(\alpha') \rangle ^{\bot }$ is $\langle 1 \rangle \oplus 8\langle -1 \rangle$, the following lemma together with the above property of $L_3'$ gives us an orthogonal basis $v_1,\dots ,v_{10}$ of $H^2(E_3';\mathbf{Z})$ such that $v_1^2=1$, $v_2^2=\dots=v_{10}^2=-1$ and $K_3'=3v_1-v_2-\dots-v_{10}$. 
\begin{lemma}[{Stipsicz-Szab\'{o} \cite[The proof of Proposition 4.3]{SS}, cf.~Wall \cite[The proof of 1.6]{W}}]\label{lem:5.5}
Let $M$ be a free $\mathbf{Z}$-module equipped with a symmetric bilinear form $\langle 1 \rangle \oplus 8\langle -1 \rangle$. If a characteristic element $K$ of $M$ satisfies $K^2=1$, then there exists an automorphism of $M$ which preserves the symmetric bilinear form on $M$ and maps $K$ to $3v_1-v_2-\dots-v_9$. Here $v_1,\dots ,v_9$ denotes an arbitrary orthogonal basis of $M$ such that $v_1^2=1$ and $v_2^2=\dots=v_9^2=-1$. 
\end{lemma}
Similarly the above lemma together with Lemma~\ref{lem:5.3} gives us an orthogonal basis $w_1,\dots ,w_{10}$ of $H^2(E(1)_{2,3};\mathbf{Z})$ such that $w_1^2=1$, $w_2^2=\dots=w_{10}^2=-1$ and $K_3=3w_1-w_2-\dots-w_{10}$. Let $\varphi ':H^2(E_3';\mathbf{Z})\to H^2(E(1)_{2,3};\mathbf{Z})$ be the isomorphism defined by $v_i\mapsto w_i\, (1\leq i\leq 10)$. The isomorphism $\varphi '$ preserves the intersection forms and the Seiberg-Witten invariants. 

Let $H\in H^2_+(R_{11};\mathbf{R})$ be the element defined in the proof of Lemma~\ref{lem:5.3}. Recall that we can view $H$ and $H'$ as positively oriented elements of $H_+^2(E(1)_{2,3};\mathbf{R})$ and $H_+^2(E_3';\mathbf{R})$, respectively. 
Note that $(-K_3)\cdot H=1$ and 
\begin{equation*}
(-K_3)\cdot \varphi'(H')=\varphi'(-K_3')\cdot \varphi'(H')=(-K_3')\cdot H'=9.
\end{equation*}
We thus have $(-K_3)\cdot H>0$ and $(-K_3)\cdot \varphi'(H')>0$. 
These two inequalities together with the lemma below show $H\cdot \varphi'(H')>0$. Hence $\varphi'$ preserves the homology orientations. 
Freedman's theorem gives us a required homeomorphism $\Phi':E(1)_{2,3}\to E_3'$ which preserves the orientations and satisfies $\Phi'^*=\varphi'$. 
\begin{lemma}\label{lem:5.6}
Let $V$ be a vector space of rank $n$ over $\mathbf{R}$ equipped with a symmetric bilinear form such that $b_2^+(V)=1$ and $b_2^-(V)=n-1$. Here $b_2^+$ and $b_2^-$ are the same notation as that in Gompf-Stipsicz~\cite{GS}. 

If elements $u, v, w\in V$ satisfy $u^2>0$, $v^2>0$, $w^2\geq 0$, $u\cdot w>0$ and $v\cdot w>0$, 
then $u\cdot v>0$.
\end{lemma}
\proof
Let $\langle u \rangle$ be the subspace spanned by $u$, and $\langle u \rangle^{\perp}$ the orthogonal complement of $\langle u \rangle$. The subspace $\langle u \rangle^{\perp}$ is negative definite, because $b_2^+(V)=1$ and $u^2>0$. Since $V=\langle u \rangle \oplus \langle u \rangle^{\perp}$, there exist elements $a,a'\in \mathbf{Z}$, $x,x'\in \langle u \rangle^{\perp}$ such that $v=au+x$ and $w=a'u+x'$. Cauchy-Schwartz inequality implies $\sqrt{x^2(x')^2}\geq x\cdot x'$. Inequalities $v^2>0$ and $w^2\geq 0$ give us $a^2u^2>-x^2$ and $(a')^2u^2\geq -(x')^2$. These three inequalities together with $v\cdot w>0$ show $aa'u^2+\lvert aa'\rvert u^2>0$. This inequality and $u\cdot w>0$ give us $a>0$. Hence $u\cdot v>0$. 

This completes the proof of Proposition~\ref{prop:5.4}. 
\end{proof}
\section{Further remarks}
\begin{figure}[htbp]
\begin{center}
\begin{minipage}{.45\linewidth}
\begin{center}
\includegraphics[scale=0.77]{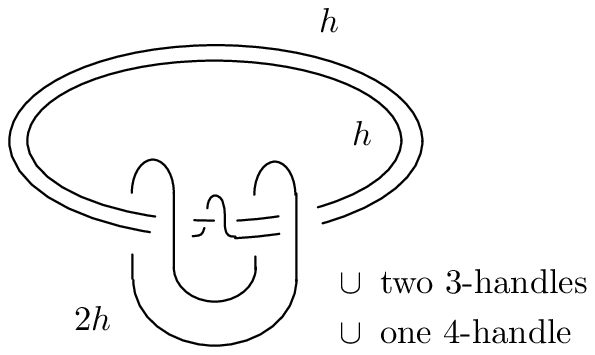}
\caption{$\mathbf{CP}^2$}
\label{4.7}
\end{center}
\end{minipage}
\begin{minipage}{.45\linewidth}
\begin{center}
\includegraphics[scale=0.77]{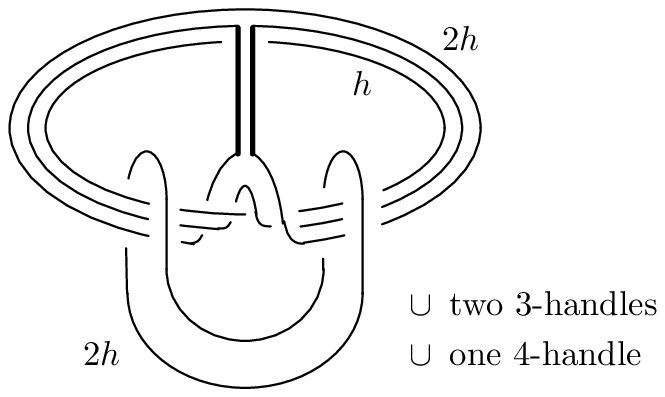}
\caption{$\mathbf{CP}^2$}
\label{4.8}
\end{center}
\end{minipage}
\end{center}
\end{figure}

\begin{figure}[htbp]
\begin{center}
\begin{minipage}{.45\linewidth}
\begin{center}
\includegraphics[scale=0.77]{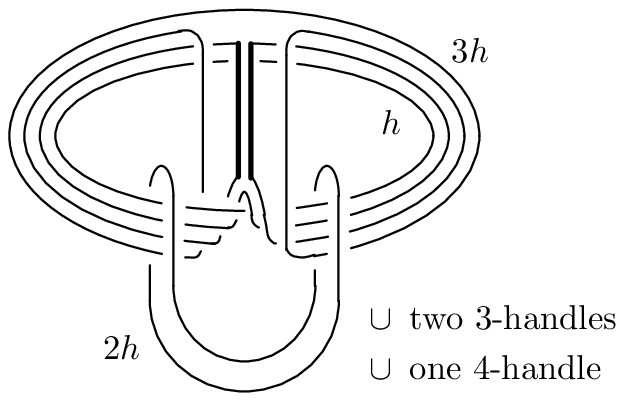}
\caption{$\mathbf{CP}^2$}
\label{4.9}
\end{center}
\end{minipage}
\begin{minipage}{.45\linewidth}
\begin{center}
\includegraphics[scale=0.68]{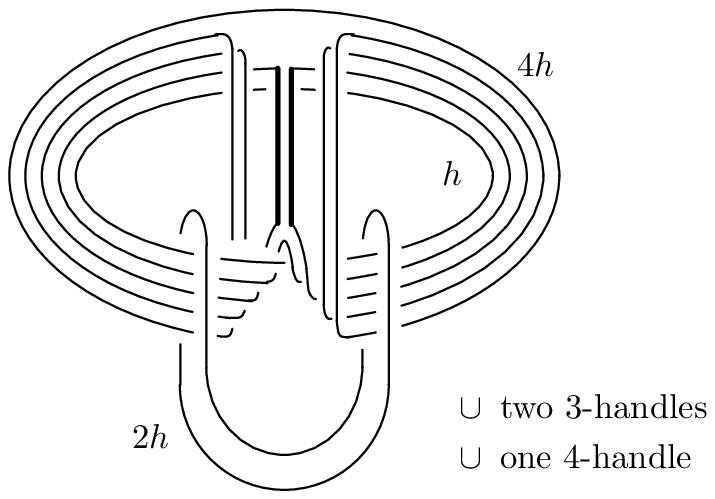}
\caption{$\mathbf{CP}^2$}
\label{4.10}
\end{center}
\end{minipage}
\end{center}
\end{figure}

We conclude this article by making some remarks. 
\begin{remark}\normalfont 
In Figure~\ref{4.8}\,$\sim $\,\ref{4.10}, we used the peculiar bands, that is, bands not in local positions to prove Lemma~\ref{lem:3.1}. Note that standard bands, that is, bands in local positions are also enough to prove Lemma~\ref{lem:3.1}. However, the peculiar bands are the key of our construction of exotic $\mathbf{CP}^2\# n\overline{\mathbf{C}\mathbf{P}^2}\,(5\leq n\leq 9)$ (see \cite{Y3}). 
In the proof of Lemma~\ref{lem:3.1}, we used two $2$-handles with framings $(2h,4h)$. Instead of these two $2$-handles, we can use two $2$-handles with framings both $(h,5h)$ and $(3h,3h)$ to prove Lemma~\ref{lem:3.1}. We can also use a $2$-handle with a framing $6h$ to construct Figure~\ref{further-construction}. In this construction, we can decrease the number of $3$-handles of $E_3$. Precisely $E_3$ admits a handle decomposition
\begin{equation*}
E_3=\text{one $0$-handle} \cup \text{eleven $2$-handles} \cup \text{one $3$-handle} \cup \text{one $4$-handle}.
\end{equation*}

We do not know if choices of the above bands and the above $2$-handles affect diffeomorphism types of $E_3$ and $E_5$. 
\begin{figure}[ht!]
\begin{center}
\includegraphics[width=2.0in]{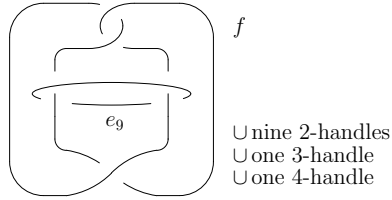}
\caption{$\mathbf{CP}^2\# 9\overline{\mathbf{C}\mathbf{P}^2}$}
\label{further-construction}
\end{center}
\end{figure}
\end{remark}
\begin{remark}\label{Yamada}
Yamada asked the author if a topologically trivial 
but smoothly non-trivial \textit{h}-cobordism between $E_q$ and $E(1)_{2,q}$ exists. 
Following the argument in Gompf-Stipsicz \cite[Example~9.2.15]{GS}, we can prove that 
such an \textit{h}-cobordism exists. Note that the same argument also shows that a topologically trivial 
but smoothly non-trivial \textit{h}-cobordism between $E(1)_{2,q}$ and itself exists.
\end{remark}
\begin{remark}
Let $X$ be a simply connected closed smooth $4$-manifold which contains a copy of $C_p$, and $X_{(p)}$ the rational blow-down of $X$ along the copy of $C_p$. Suppose that $X_{(p)}$ is simply connected. Do the following two conditions, $X$ and the homomorphism $H_2(C_p;\mathbf{Z})\to H_2(X;\mathbf{Z})$ induced by the copy of $C_p$, suffice to determine the (small perturbation) Seiberg-Witten invariant of $X_{(p)}$? 

The proofs of Proposition~\ref{prop:5.2} and Proposition~\ref{prop:5.4} give an affirmative answer to this question in some cases. In a forthcoming paper, we will give a more general result for this question. 
\end{remark}

We here give a proof of the following proposition referred in the introduction of this article.
\begin{proposition}\label{S^4}
If a smooth $4$-manifold is homeomorphic to $\mathbf{S}^4$ {\normalfont (}resp.~$\mathbf{CP}^2${\normalfont )} and admits neither $1$- nor $3$-handles in a handle decomposition, then the $4$-manifold is diffeomorphic to $\mathbf{S}^4$ {\normalfont (}resp.~$\mathbf{CP}^2${\normalfont )}.
\end{proposition}
\begin{proof}Note that if a simply connected closed smooth $4$-manifold has neither $1$- nor $3$-handles in a handle decomposition, then the number of $2$-handles appeared in the handle decomposition is equal to the rank of the second homology group of the $4$-manifold. 

Suppose that a smooth $4$-manifold is homeomorphic to $\mathbf{S}^4$ and has neither $1$- nor $3$-handles in a handle decomposition. 
Then this handle body consists of a $0$-handle and a $4$-handle. Since attaching a $4$-handle is unique (see Gompf-Stipsicz \cite{GS}), the $4$-manifold is diffeomorphic to $\mathbf{S}^4$. 

Suppose that a smooth $4$-manifold is homeomorphic to $\mathbf{CP}^2$ and has neither $1$- nor $3$-handles in a handle decomposition. Then this handle body consists of a $0$-handle, a $2$-handle and a $4$-handle. Thus the attaching circle of the $2$-handle produces $\mathbf{S}^3$ by a Dehn surgery with coefficient $+1$. Since such a knot is unknot (see Gordon-Luecke \cite{GL}), the $4$-manifold is diffeomorphic to $\mathbf{CP}^2$. 
\end{proof}
Contrary to the above proposition, many simply connected closed topological $4$-manifolds are known to admit at least two different smooth structures without $1$- and $3$-handles (cf.~Gompf-Stipsicz \cite{GS}). As far as the author knows, $\mathbf{S}^4$ and $\mathbf{CP}^2$ are the only known exceptions. Thus the following problem is natural. 
\begin{problem}\normalfont 
Which simply connected closed topological $4$-manifold has a unique smooth structure without $1$- and $3$-handles?
\end{problem}
Finally we refer to further constructions. 
\begin{remark}\normalfont 
This article is based on the author's announcement \cite{Y1}. In \cite{Y2}, we will give the rest of examples announced in \cite{Y1}. In addition to these examples, we will construct a smooth $4$-manifold which has the same Seiberg-Witten invariant as $E(1)_{2,3}$ and admits no $1$-handles as follows: We `naturally' construct Figure~\ref{further-construction} and perform a logarithmic transformation of multiplicity $3$ in the cusp neighborhood.

 In \cite{Y3}, we will construct examples of exotic $\mathbf{CP}^2\# n\overline{\mathbf{C}\mathbf{P}^2}\,(5\leq n\leq 9)$ by using rational blow-downs and Kirby calculus. We also prove that our examples admit a handle decomposition without $1$- and $3$-handles in the case $7\leq n\leq 9$. 
\end{remark}



\end{document}